\RequirePackage{ifpdf}
\ifpdf 
\documentclass[pdftex]{sigma}
\else
\documentclass{sigma}
\fi

\numberwithin{equation}{section}
\numberwithin{theorem}{section}
\numberwithin{proposition}{section}
\numberwithin{lemma}{section}
\numberwithin{corollary}{section}
\numberwithin{definition}{section}
\numberwithin{example}{section}
\numberwithin{remark}{section}
\numberwithin{note}{section}

\newcommand{\comment}[1]{}
\newcommand{\dsum}[2]{{\displaystyle\sum_{#1}^{#2}}}
\newcommand{\dprod}[2]{{\displaystyle\prod_{#1}^{#2}}}

\newcommand{\pr}[1]{\{\,#1\,\}}
\newcommand{\br}[1]{\langle{#1}\rangle}
\newcommand{\sgm}[1]{[#1]}
\newcommand{\sbinom}[2]{\left[\begin{matrix}{#1}\\{#2}\end{matrix}\right]}
\newcommand{\cA}{\mathcal{A}}
\newcommand{\cE}{\mathcal{E}}
\newcommand{\cR}{\mathcal{R}}
\newcommand{\cG}{\mathcal{G}}

\begin{document}

\allowdisplaybreaks

\renewcommand{\thefootnote}{$\star$}

\renewcommand{\PaperNumber}{054}

\FirstPageHeading

\ShortArticleName{Kernel Functions for Dif\/ference Operators of Ruijsenaars Type}

\ArticleName{Kernel Functions for Dif\/ference Operators\\ of
Ruijsenaars Type
and Their Applications\footnote{This paper is a contribution to the Proceedings of the Workshop ``Elliptic Integrable Systems, Isomonodromy Problems, and Hypergeometric Functions'' (July 21--25, 2008, MPIM, Bonn, Germany). The full collection
is available at
\href{http://www.emis.de/journals/SIGMA/Elliptic-Integrable-Systems.html}{http://www.emis.de/journals/SIGMA/Elliptic-Integrable-Systems.html}}}

\Author{Yasushi KOMORI~$^\dag$, Masatoshi NOUMI~$^\ddag$ and Jun'ichi SHIRAISHI~$^\S$}

\AuthorNameForHeading{Y.~Komori, M.~Noumi and J.~Shiraishi}

\Address{$^\dag$~Graduate School of Mathematics, Nagoya University,
Chikusa-Ku, Nagoya 464-8602, Japan}
\EmailD{\href{mailto:komori@math.nagoya-u.ac.jp}{komori@math.nagoya-u.ac.jp}}

\Address{$^\ddag$~Department of Mathematics, Kobe University,
Rokko, Kobe 657-8501, Japan}
\EmailD{\href{mailto:noumi@math.kobe-u.ac.jp}{noumi@math.kobe-u.ac.jp}}

\Address{$^\S$~Graduate School of Mathematical Sciences, University of Tokyo,\\
\hphantom{$^\S$}~Komaba, Tokyo 153-8914, Japan}
\EmailD{\href{mailto:shiraish@ms.u-tokyo.ac.jp}{shiraish@ms.u-tokyo.ac.jp}}

\ArticleDates{Received December 01, 2008, in f\/inal form April 30,
2009; Published online May 12, 2009}

\Abstract{A unif\/ied approach is given to kernel functions which intertwine Ruijsenaars dif\/ference operators of type $A$ and of type $BC$.
As an application of the trigonometric cases, new explicit formulas for Koornwinder
polynomials attached to single columns and single rows are derived.}


\Classification{81R12; 33D67}

\tableofcontents

\renewcommand{\thefootnote}{\arabic{footnote}}
\setcounter{footnote}{0}

\section{Introduction}\label{section1}

Let $\Phi(x;y)$ be a meromorphic function on
$\mathbb{C}^m\times\mathbb{C}^n$ in the variables
$x=(x_1,\ldots,x_m)$ and $y=(y_1,\ldots,y_n)$,
and consider two operators $\mathcal{A}_x$, $\mathcal{B}_y$
which act on meromorphic functions in $x$ and $y$, respectively.
We say that $\Phi(x;y)$ is a {\em kernel function} for the
pair $(\mathcal{A}_x,\mathcal{B}_y)$ if it satisf\/ies a
functional equation of the form
\begin{gather*}
\mathcal{A}_x \Phi(x;y)=\mathcal{B}_y\Phi(x;y).
\end{gather*}
In the theory of Jack and Macdonald polynomials \cite{Ma1995},
certain explicit kernel functions
play crucial roles in eigenfunction expansions
and integral representations.
Recently, kernel functions in this sense have been studied
systematically by Langmann \cite{Lm2004, Lm2006} in the analysis of
eigenfunctions for the elliptic quantum integrable systems of
Calogero--Moser type, and by Ruijsenaars \cite{Ru2005, Ru2006, Ru2009} for the
relativistic elliptic quantum integrable systems of
Ruijsenaars--Schneider type.

In this paper we investigate two kinds of kernel functions,
{\em of Cauchy type} and {\em of dual Cauchy type},
which intertwine pairs of Ruijsenaars dif\/ference operators.
In the cases of elliptic dif\/ference operators,
kernel functions of Cauchy type for the ($A_{n-1},A_{n-1}$) and
($BC_n,BC_n$) cases were found by
Ruijsenaars~\cite{Ru2006, Ru2009}.
Extending his result,
we present kernel functions of Cauchy type,
as well as those of dual Cauchy type,
for the ($BC_m,BC_n$) cases
with arbitrary $m$, $n$
(under certain balancing conditions on the parameters in the elliptic cases).
For the trigonometric dif\/ference operators
of type $A$, kernel functions both of Cauchy and
dual Cauchy types were already discussed by
Macdonald~\cite{Ma1995}.
Kernel functions of dual Cauchy type for the trigonometric $BC_n$
cases are due to Mimachi~\cite{Mi2001}.
In this paper we develop a unif\/ied approach to kernel functions
for Ruijsenaars operators of type $A$ and of type $BC$,
with rational, trigonometric and elliptic coef\/f\/icients,
so as to cover all these known examples in the dif\/ference cases.
We expect that our framework could be ef\/fectively applied
to the study of eigenfunctions for dif\/ference operators of Ruijsenaars
type.

As such an application of kernel functions in the trigonometric
$BC$ cases, we derive new explicit formulas for
Koornwinder polynomials attached to single columns and single rows.
This provides with
a direct construction of those special cases of
the binomial expansion formula for the Koornwinder polynomials
as studied by Okounkov~\cite{O1998}
and Rains~\cite{Ra2005}.
We also remark that,
regarding explicit formulas for Macdonald polynomials
attached to single rows of type~$B$, $C$, $D$,
some conjectures have been proposed by Lassalle~\cite{L2004}.
The relationship between his conjectures and our kernel functions
will be discussed in a separate paper.

Our main results on the kernel functions for Ruijsenaars dif\/ference
operators of type $A$ and of type $BC$ will be formulated
in Section~\ref{section2}. In Section~\ref{section3} we give a unif\/ied proof
for them on the basis of two key identities.
After giving remarks in Section~\ref{section4}
on the passage to the
$q$-dif\/ference operators of Macdonald and Koornwinder,
as an application of our approach
we present in Section~\ref{section5} explicit formulas for Koornwinder polynomials
attached to single columns and single rows.
We also include three sections in Appendix.
We will give some remarks in Appendix~\ref{appendixA} on higher order dif\/ference operators,
and in Appendix~\ref{appendixB} make an explicit
comparison of our kernel functions in the elliptic
cases with those constructed by Ruijsenaars \cite{Ru2009}.
In Appendix~\ref{appendixC},
we will give a proof of the fact that certain Laurent polynomials,
which appear in our explicit formulas for Koornwinder polynomials
attached to single columns and single rows,
are special cases of the $BC_m$ interpolation polynomials
of Okounkov~\cite{O1998}.

\section{Kernel functions for Ruijsenaars operators}\label{section2}

\subsection{Variations of the gamma function}\label{section2.1}

In order to specify the class of
Ruijsenaars operators which we shall discuss below,
by the symbol~$\sgm{u}$ we denote a nonzero entire function
in one variable $u$, satisfying the following
{\em Riemann relation}:
\begin{gather}
\sgm{x+u}\sgm{x-u}\sgm{y+v}\sgm{y-v}
-
\sgm{x+v}\sgm{x-v}\sgm{y+u}\sgm{y-u}\nonumber\\
\qquad{} =
\sgm{x+y}\sgm{x-y}\sgm{u+v}\sgm{u-v}\label{eq:Riemann}
\end{gather}
for any $x,y,u,v\in\mathbb{C}$.
Under this condition, it is known that
$\sgm{u}$ must be an odd function
($\sgm{-u}=-\sgm{u}$, and hence $\sgm{0}=0$), and
the set $\Omega$ of all zeros of~$\sgm{u}$
forms an additive subgroup of~$\mathbb{C}$.
Furthermore such functions are classif\/ied
into the three categories,
{\em rational}, {\em trigonometric} and {\em elliptic},
according to the rank of the additive subgroup $\Omega\subset\mathbb{C}$
of all zeros of $\sgm{u}$. In fact, up to multiplication by nonzero constants,
$\sgm{u}$ coincides with one of the following three types of functions:
\begin{alignat*}{4}
& (0)\ \mbox{rational case:} \qquad && e\left(a u^2\right)u \qquad &&(\Omega=0),& \\
& (1)\ \mbox{trigonometric case:} \qquad& & e\left(a u^2\right)\sin(\pi u/\omega_1)\qquad && (\Omega=\mathbb{Z}\omega_1),&\\
& (2)\ \mbox{elliptic case:}  \qquad && e\left(au^2\right)\sigma(u;\Omega)\qquad &&(\Omega=\mathbb{Z}\omega_1\oplus\mathbb{Z}\omega_2),&
\end{alignat*}
where $a\in\mathbb{C}$ and
$e(u)=\exp(2\pi\sqrt{-1}u)$.  Also,
$\sigma(x;\Omega)$ denotes the Weierstrass sigma function
\begin{gather*}
\sigma(u;\Omega)=
u\dprod{\omega\in\Omega;\, \omega\ne 0}{}
\left(1-\tfrac{u}{\omega}\right) e^{\frac{u}{\omega}+\frac{u^2}{2\omega^2}}
\end{gather*}
associated with the period lattice $\Omega=\mathbb{Z}\omega_1\oplus
\mathbb{Z}\omega_2$, generated by $\omega_1$, $\omega_2$
which are linearly independent over $\mathbb{R}$.

We start with some remarks on gamma functions associated with
the function $\sgm{u}$.
Fixing a nonzero scaling constant $\delta\in\mathbb{C}$,
suppose that a nonzero meromorphic function
$G(u|\delta)$ on $\mathbb{C}$ satisf\/ies the dif\/ference equation
\begin{gather}\label{eq:ellgam}
G(u+\delta|\delta)=\sgm{u} G(u|\delta)\qquad(u\in\mathbb{C}).
\end{gather}
Such a function $G(u|\delta)$, determined up to multiplication
by $\delta$-periodic functions, will be called a
{\em gamma function} for $\sgm{u}$.
We give typical examples of gamma functions in this sense
for the rational, trigonometric and elliptic cases.

{\bf (0) Rational case:}
For any $\delta\in\mathbb{C}$ ($\delta\ne 0$),
the meromorphic function{\samepage
\begin{gather*}
G(u|\delta)=\delta^{u/\delta} \Gamma(u/\delta),
\end{gather*}
def\/ined with the Euler gamma function $\Gamma(u)$, is a
gamma function for $\sgm{u}=u$.}

{\bf (1) Trigonometric case:}
We set $z{=}e(u/\omega_1)$ and $q{=}e(\delta/\omega_1)$,
and suppose that \mbox{$\mbox{Im}(\delta/\omega_1){>}0$} so that $|q|<1$.
We now consider the function
\begin{gather*}
\sgm{u}=z^{{\frac 12}}-z^{-{\frac 12}}=2\sqrt{-1}\sin(\pi u/\omega_1).
\end{gather*}
For this $\sgm{u}$ two meromorphic functions
\begin{gather}\label{eq:tgamma}
G_{-}(u|\delta)=
\dfrac{e\big(-\tfrac{\delta}{2\omega_1}\tbinom{u/\delta}{2}\big)}{(z;q)_\infty},
\qquad
G_{+}(u|\delta)=
e\big(\tfrac{\delta}{2\omega_1}\tbinom{u/\delta}{2}\big) (q/z;q)_\infty,
\end{gather}
where $(z;q)_\infty=\prod\limits_{i=0}^{\infty}(1-q^iz)$,
satisfy the dif\/ference equations
\begin{gather*}
G_{-}(u+\delta|\delta)=-\sgm{u}  G_{-}(u|\delta),
\qquad
G_{+}(u+\delta|\delta)=\sgm{u}  G_{+}(u|\delta),
\end{gather*}
respectively.
Namely, for $\epsilon=\pm$,
$G_\epsilon(u|\delta)$ is a gamma function for $\epsilon \sgm{u}$.
(Note that the quadratic function $\tbinom{u}{2}={\frac 12} u(u-1)$
satisf\/ies $\tbinom{u+1}{2}=\tbinom{u}{2}+u$.)
Another set of gamma functions for $\epsilon\sgm{u}$ of this case
is given by
\begin{gather*}
G_{+}(u|\delta)
=\dfrac{\big(q^{{\frac 12}}z^{-{\frac 12}};q^{{\frac 12}}\big)_\infty}
{\big(-z^{{\frac 12}};q^{{\frac 12}}\big)_\infty},
\qquad
G_{-}(u|\delta)
=\dfrac{\big(-q^{{\frac 12}}z^{-{\frac 12}};q^{{\frac 12}}\big)_\infty}
{\big(z^{{\frac 12}};q^{{\frac 12}}\big)_\infty}
\end{gather*}
by means of the inf\/inite products with base $q^{{\frac 12}}$.

{\bf (2) Elliptic case:}
Let $p$, $q$ be nonzero complex numbers with $|p|<1$, $|q|<1$.  Then
the {\em Ruijsenaars elliptic gamma function}
\begin{gather*}
\Gamma(z;p,q)=\dfrac{(pq/z;p,q)_\infty}{(z;p,q)_\infty},
\qquad
(z;p,q)_\infty=\dprod{i,j=0}{\infty}\left(1-p^iq^jz\right)
\end{gather*}
satisf\/ies the $q$-dif\/ference equation
\begin{gather*}
\Gamma(qz;p,q)=\theta(z;p) \Gamma(z;p,q),
\qquad\theta(z;p)=(z;p)_\infty(p/z;p)_\infty.
\end{gather*}
Note also $\Gamma(pq/z;p,q)=\Gamma(z;p,q)^{-1}$.
We set $p=e(\omega_2/\omega_1)$,
$q=e(\delta/\omega_1)$ and $z=e(u/\omega_1)$,
and suppose that $\mbox{Im}(\omega_2/\omega_1)>0$,
$\mbox{Im}(\delta/\omega_1)>0$ so that $|p|<1$, $|q|<1$.
Instead of $\theta(z;p)$ above,
we consider the function
\begin{gather*}
\sgm{u}
=-z^{-{\frac 12}}\theta(z;p)
=2\sqrt{-1}\sin(\pi u/\omega_1)
(pz;p)_\infty (p/z;p)_\infty
\end{gather*}
which is a constant multiple of the odd Jacobi theta
function with modulus $\omega_2/\omega_1$,
so that $\sgm{u}$ should satisfy the Riemann relation.
Then the two meromorphic functions
\begin{gather}
G_{-}(u|\delta)=
e\big(-\tfrac{\delta}{2\omega_1}\tbinom{u/\delta}{2}\big)
\Gamma(z;p,q)
=
e\big(-\tfrac{\delta}{2\omega_1}\tbinom{u/\delta}{2}\big)
\dfrac{(pq/z;p,q)_\infty}{(z;p,q)_\infty},\nonumber\\
G_{+}(u|\delta)=
e\big(\tfrac{\delta}{2\omega_1}\tbinom{u/\delta}{2}\big)
\Gamma(pz;p,q)
=
e\big(\tfrac{\delta}{2\omega_1}\tbinom{u/\delta}{2}\big)
\dfrac{(q/z;p,q)_\infty}{(pz;p,q)_\infty}.\label{eq:egammapm}
\end{gather}
satisfy the dif\/ference equations
\begin{gather*}
G_{-}(u+\delta|\delta)=-\sgm{u}  G_{-}(u|\delta),
\qquad
G_{+}(u+\delta|\delta)=\sgm{u}  G_{+}(u|\delta),
\end{gather*}
respectively.
Another set of gamma functions for $\epsilon\sgm{u}$ is given by
\begin{gather*}
G_{\epsilon}(u|\delta)=\theta\big(\epsilon z^{{\frac 12}};q^{{\frac 12}}\big)\Gamma(z;p,q)
=
\Gamma\big(\epsilon p^{{\frac 12}}z^{{\frac 12}};p^{{\frac 12}},q^{{\frac 12}}\big)
\Gamma\big(-\epsilon z^{{\frac 12}};p^{{\frac 12}},q^{{\frac 12}}\big)\qquad(\epsilon=\pm).
\end{gather*}
In the limit as $p\to 0$, these examples
recover the previous ones in the trigonometric case.

We remark that, if $G(u|\delta)$ is a gamma function for $\sgm{u}$, then
$G(\delta-u|\delta)^{-1}$ is a gamma function for $-\sgm{u}$.
Also, when we transform $\sgm{u}$ to $\sgm{u}'=c e(au^2)\sgm{u}$,
a gamma function for $\sgm{u}'$ is obtained for instance
as $G'(u|\delta)=c^{u/\delta}e(a p(u))G(u|\delta)$,
where $p(u)=\tfrac{1}{3\delta}u^3-\tfrac{1}{2}u^2+\tfrac{\delta}{6}u$.

\subsection[Kernel functions of type $A$]{Kernel functions of type $\boldsymbol{A}$}\label{section2.2}

The elliptic dif\/ference operator which we will discuss below
was introduced by Ruijsenaars~\cite{Ru1987} together with the
commuting family of higher order dif\/ference operators.
In order to deal with rational, trigonometric and elliptic cases
in a unif\/ied manner,
we formulate our results for this class of (f\/irst order) dif\/ference
operators in terms of an arbitrary function $\sgm{u}$ satisfying
the Riemann relation.
(As to the commuting family of higher order dif\/ference operators,
we will give some remarks later in Appendix~\ref{appendixA}.)

Fix a nonzero entire function $\sgm{u}$ satisfying
the Riemann relation \eqref{eq:Riemann}.
For type $A$, we consider the dif\/ference operator
\begin{gather*}
D_{x}^{(\delta,\kappa)}
=\dsum{i=1}{m}  A_i(x;\kappa)  T_{x_i}^{\delta},
\qquad
A_i(x;\kappa)=
\dprod{1\le j\le m;\,j\ne i}{}
\dfrac{\sgm{x_i-x_j+\kappa}}{\sgm{x_i-x_j}}
\end{gather*}
in $m$ variables $x=(x_1,\ldots,x_m)$,
where $\delta,\kappa\in\mathbb{C}$ are complex parameters $\notin\Omega$,
and
$T_{x_i}^{\delta}$ stands for the $\delta$-shift operator
\begin{gather*}
T_{x_i}^{\delta}f(x_1,\ldots,x_m)=
f(x_1,\ldots,x_i+\delta,\ldots,x_m)\qquad(i=1,\ldots,m).
\end{gather*}
Note that this operator remains invariant if one replaces the
function $\sgm{u}$ with its multiple by any nonzero constant.

By taking a gamma function $G(u|\delta)$ for (any constant multiple of) $\sgm{u}$
as in \eqref{eq:ellgam},
we def\/ine a function $\Phi_{A}(x;y|\delta,\kappa)$ by
\begin{gather}\label{eq:phifun}
\Phi_{A}(x;y|\delta,\kappa)
=\dprod{j=1}{m}\dprod{l=1}{n} \dfrac{G(x_j+y_l+v-\kappa|\delta)}
{G(x_j+y_l+v|\delta)}
\end{gather}
with an extra parameter $v$.
We also consider the function
\begin{gather}\label{eq:psifun}
\Psi_{ A}(x;y)=\dprod{j=1}{m}\dprod{l=1}{n}
 \sgm{x_j-y_l+v}.
\end{gather}
These two functions
$\Phi_{ A}(x;y|\delta,\kappa)$
and
$\Psi_{ A}(x;y)$
are kernel functions {\em of Cauchy type} and {\em of dual
Cauchy type} for this case, respectively.

\begin{theorem}\label{thm:AE}
Let $\sgm{u}$ be any nonzero entire function satisfying
the Riemann relation.
\begin{enumerate}\itemsep=0pt
\item[$(1)$] If $m=n$, then
the function $\Phi_A(x;y|\delta,\kappa)$
defined as \eqref{eq:phifun} satisfies the functional equation
\begin{gather*}
D_x^{(\delta,\kappa)}\Phi_{A}(x;y|\delta,\kappa)=
D_y^{(\delta,\kappa)}\Phi_{A}(x;y|\delta,\kappa).
\end{gather*}
\item[$(2)$]
Under the balancing condition $m\kappa+n\delta=0$,
the function $\Psi_A(x;y)$ defined as \eqref{eq:psifun} satisfies
the functional equation
\begin{gather*}
\sgm{\kappa} D_x^{(\delta,\kappa)}\Psi_{A}(x;y)+
\sgm{\delta} D_y^{(\kappa,\delta)}\Psi_{A}(x;y)=0.
\end{gather*}
\end{enumerate}
\end{theorem}

Statement (1) of Theorem \ref{thm:AE} is due to Ruijsenaars \cite{Ru2006,Ru2009}.
(See Appendix~\ref{appendixB.1} for an explicit comparison between
$\Phi_{A}(x;y)$ and Ruijsenaars' kernel function of~\cite{Ru2009}.)
In the scope of the present paper,
the {\em balancing conditions}
($m=n$ in~(1), and $m\kappa+n\delta=0$ in~(2))
seem to be essential in the elliptic cases.
In the context of elliptic dif\/ferential operators of Calogero--Moser type,
however,
Langmann~\cite{Lm2006} has found
a natural generalization of the kernel identities of Cauchy type,
which include the dif\/ferentiation with
respect to the elliptic modulus, to arbitrary pair $(m,n)$.
It would be a intriguing problem to f\/ind a generalization
of this direction for elliptic dif\/ference operators of Ruijsenaars type.

In trigonometric and rational cases,
these functions
$\Phi_{A}(x;y|\delta,\kappa)$ and $\Psi_{A}(x;y)$ satisfy more general
functional equations without balancing conditions.
\begin{theorem}\label{thm:AT}
Suppose that $\sgm{u}$ is a constant multiple
of $\sin(\pi u/\omega_1)$ or $u$.
\begin{enumerate}\itemsep=0pt
\item[$(1)$] For arbitrary $m$ and $n$,
the function $\Phi_{A}(x;y|\delta,\kappa)$
satisfies the functional equation
\begin{gather*}
D_x^{(\delta,\kappa)}\Phi_{A}(x;y|\delta,\kappa)-
D_y^{(\delta,\kappa)}\Phi_{A}(x;y|\delta,\kappa)
=\dfrac{\sgm{(m-n)\kappa}}{\sgm{\kappa}}
\Phi_{A}(x;y|\delta,\kappa).
\end{gather*}
\item[$(2)$]
The function $\Psi_{A}(x;y)$ satisfies
the functional equation
\begin{gather*}
\sgm{\kappa} D_x^{(\delta,\kappa)}\Psi_{A}(x;y)+
\sgm{\delta} D_y^{(\kappa,\delta)}\Psi_{A}(x;y)
=\sgm{m\kappa+n\delta}\Psi_{A}(x;y).
\end{gather*}
\end{enumerate}
\end{theorem}
These results for the trigonometric (and rational) cases are
essentially contained in the discussion of Macdonald \cite{Ma1995}.

A unif\/ied proof of Theorems \ref{thm:AE} and \ref{thm:AT} will
be given in Section~\ref{section3}.
We will also explain in Section~\ref{section4} how Theorem~\ref{thm:AT} is
related with the theory of Macdonald polynomials.

\subsection[Kernel functions of type $BC$]{Kernel functions of type $\boldsymbol{BC}$}\label{section2.3}

The (f\/irst-order) elliptic dif\/ference operator of type $BC$
was f\/irst proposed by van Diejen \cite{vD1994}.
It is also known by Komori--Hikami \cite{KH1997} that it
admits a commuting family of higher order dif\/ference operators.
In the following we use the expression of the
f\/irst-order dif\/ference operator due to \cite{KH1997},
with modif\/ication in terms of $\sgm{u}$.
(In Appendix~\ref{appendixB}, we will give some remarks on the comparison
of our dif\/ference operator with other expressions in the
literature.)

For type $BC$, we consider dif\/ference operators of the form
\begin{gather}\label{eq:defE}
E_{x}^{(\mu|\delta,\kappa)}
=
\dsum{i=1}{m}
A_{i}^{+}(x;\mu|\delta,\kappa) T_{x_i}^{\delta}
+
\dsum{i=1}{m}
A_{i}^{-}(x;\mu|\delta,\kappa) T_{x_i}^{-\delta}
+A^{0}(x;\mu|\delta,\kappa)
\end{gather}
including $2\rho$ parameters $\mu=(\mu_1,\ldots,\mu_{2\rho})$
besides $(\delta,\kappa)$,
where $\rho=1,\, 2$ or $4$ according as $\mbox{\rm rank}\,\Omega=0,\, 1$,
or $2$.
In the trigonometric and rational cases of type $BC$,
we assume that the function $\sgm{u}$ does {\em not} contain exponential
factors. Namely, we assume that $\sgm{u}$ is a constant multiple of
one of the functions
\begin{gather}\label{eq:sgmBC}
u\ \ (\rho=1),\qquad
\sin(\pi u/\omega_1)\ \ (\rho=2),\qquad
e\left(au^2\right) \sigma(u;\Omega)\ \ (\rho=4).
\end{gather}
In each case we def\/ine $\omega_1,\ldots,\omega_\rho\in\Omega$ as
\begin{alignat*}{4}
& \mbox{(0)\ rational case:} \qquad && \Omega=0, \qquad && \omega_1=0,&
\\
& \mbox{(1)\ trigonometric case:} \qquad && \Omega=\mathbb{Z}\omega_1,\qquad && \omega_2=0,&
\\
& \mbox{(2)\ elliptic case:} \qquad && \Omega=\mathbb{Z}\omega_1\oplus
\mathbb{Z}\omega_2, \qquad && \omega_3=-\omega_1-\omega_2,\quad \omega_4=0.&
\end{alignat*}
Then the quasi-periodicity of $\sgm{u}$ is described as
\begin{gather*}
\sgm{u+\omega_r}=\epsilon_r  e\big(\eta_r\big(u+{\tfrac 12}\omega_r\big)\big)\sgm{u}
\qquad (r=1,\ldots,\rho)
\end{gather*}
for some $\eta_r\in\mathbb{C}$ and $\epsilon_r=\pm1$.
In the trigonometric and rational cases (without exponential factors),
one can simply take $\eta_r=0$ ($r=1,\ldots,\rho$).
Note also that $\sgm{u}$ admits the duplication formula of the form
\begin{gather}\label{eq:dup}
\sgm{2u}=2\sgm{u}
\dprod{s=1}{\rho-1}
\dfrac{\sgm{u-{\frac 12}\omega_s}}
{\sgm{-{\frac 12}\omega_s}},
\qquad
\dfrac{\sgm{2u+c}}{\sgm{2u}}=\dprod{s=1}{\rho}
\dfrac{\sgm{u+{\frac 12}(c-\omega_s)}}{\sgm{u-{\frac 12}\omega_s}}.
\end{gather}
(In the trigonometric and rational cases,
these formulas fail for $\sgm{u}$ containing
nontrivial exponential factors $e(au^2)$.)
In relation to the parameters $\mu=(\mu_1,\ldots,\mu_{2\rho})$,
we introduce
\begin{gather*}
c^{(\mu|\delta,\kappa)}=\sum_{s=1}^{2\rho}\mu_s
-\tfrac{\rho}{2}(\delta+\kappa)
+\sum_{s=1}^{\rho}\omega_s.
\end{gather*}
This parameter is related to quasi-periodicity of
the coef\/f\/icients of the Ruijsenaars operator, and
plays a crucial role in various places of our argument.
Note that the last term $\sum\limits_{s=1}^{\rho}\omega_s$ is nontrivial
only in the trigonometric case:
it is $\omega_1$ when $\rho=2$, and $0$ when $\rho=1,\,4$.

With these data,
we def\/ine the coef\/f\/icients of the Ruijsenaars operator
$E_{x}^{(\mu|\delta,\kappa)}$ of type $BC$ as follows:
\begin{gather}
A_i^{+}(x;\mu|\delta,\kappa)=
\dfrac{\prod\limits_{s=1}^{2\rho} \sgm{x_i+\mu_s}
}{
\prod\limits_{s=1}^{\rho}
\sgm{x_i - {\frac 12}\omega_s}
\sgm{x_i + {\frac 12}(\delta - \omega_s)}
}
\dprod{1\le j\le m;\,j\ne i}{}
\dfrac{\sgm{x_i \pm x_j + \kappa}}
{\sgm{x_i \pm x_j}},
\nonumber
\\
A_i^{-}(x;\mu|\delta,\kappa)=A_i^{+}(-x;\mu|\delta,\kappa) \qquad(i=1,\ldots,m), \label{eq:defE1}
\end{gather}
with an abbreviated notation $\sgm{x\pm y}=\sgm{x+y}\sgm{x-y}$ of products,
and
\begin{gather}
A^{0}(x;\mu|\delta,\kappa)=\dsum{r=1}{\rho}A_r^{0}(x;\mu|\delta,\kappa),
\nonumber\\
A^{0}_r(x;\mu|\delta,\kappa) =
\dfrac{e(-(m\kappa+{\frac 12} c^{(\mu|\delta,\kappa)})\eta_r)
\prod\limits_{s=1}^{2\rho} \sgm{{\frac 12}(\omega_r - \delta) + \mu_s}}
{\sgm{\kappa}\
\prod\limits_{s\ne r}
\sgm{{\frac 12}(\omega_r - \omega_s)}
\prod\limits_{s=1}^{\rho}
\sgm{{\frac 12}(\omega_r - \omega_s + \kappa - \delta)}}
\dprod{j=1}{m}
\dfrac{\sgm{{\frac 12}(\omega_r - \delta)\pm x_j + \kappa}}
{\sgm{{\frac 12}(\omega_r - \delta)\pm x_j}}.\label{eq:defE2}
\end{gather}
This operator $E_{x}^{(\mu|\delta,\kappa)}$ is one of the
expressions for the f\/irst-order Ruijsenaars operator of
type $BC_m$ due to Komori--Hikami \cite[(4.21)]{KH1997}.
We remark that this dif\/ference operator $E_{x}^{(\mu|\delta,\kappa)}$
has symmetry
\begin{gather}\label{eq:symE}
E_{x}^{(\mu|-\delta,-\kappa)}=E_{x}^{(-\mu|\delta,\kappa)}
\end{gather}
with respect to the sign change of parameters.
(This property can be verif\/ied directly by using the quasi-periodicity
of $\sgm{u}$ and the fact that $\prod\limits_{s=1}^{\rho}\epsilon_s=1,-1,-1$
according as $\rho=1,2,4$.)
Note also that it remains invariant if one replace the
function $\sgm{u}$ with its multiple by a nonzero constant.
By the duplication formula \eqref{eq:dup}, one can also rewrite the
operator $E_{x}^{(\mu|\delta,\kappa)}$ in the form{\samepage
\begin{gather}
\tfrac{1}{4}\dprod{s=1}{\rho-1}\sgm{{\tfrac 12}\omega_s}^2E_{x}^{(\mu|\delta,\kappa)}
=
\dsum{i=1}{m}
\dfrac{\prod\limits_{s=1}^{2\rho}\sgm{x_i+\mu_s}}{\sgm{2x_i}
\sgm{2x_i+\delta}}
\dprod{1\le j\le m;\,j\ne i}{}
\dfrac{\sgm{x_i\pm x_j+\kappa}}{\sgm{x_i\pm x_j}}
 T_{x_i}^\delta
\nonumber\\
\hphantom{\dfrac{1}{4}\dprod{s=1}{\rho-1}\sgm{{\frac 12}\omega_s}^2E_{x}^{(\mu|\delta,\kappa)}=}{}
+
\dsum{i=1}{m}
\dfrac{\prod\limits_{s=1}^{2\rho}\sgm{-x_i+\mu_s}}{\sgm{-2x_i}\sgm{-2x_i+\delta}}
\dprod{1\le j\le m;\,j\ne i}{}
\dfrac{\sgm{-x_i\pm x_j+\kappa}}{\sgm{-x_i\pm x_j}}
T_{x_i}^{-\delta}
\nonumber\\
\hphantom{\dfrac{1}{4}\dprod{s=1}{\rho-1}\sgm{{\frac 12}\omega_s}^2E_{x}^{(\mu|\delta,\kappa)}=}{}+
\dsum{r=1}{\rho}
K_r
\dfrac{\prod\limits_{s=1}^{2\rho}\sgm{{\frac 12}(\omega_r - \delta) + \mu_s}}
{2\sgm{\kappa}\sgm{\kappa-\delta}}
\dprod{j=1}{m}
\dfrac{\sgm{{\frac 12}(\omega_r - \delta)\pm x_j+\kappa}}{\sgm{{\frac 12}(\omega_r - \delta)\pm x_j}},
\label{eq:defEE}
\end{gather}
where $K_r=e\big(-(\omega_r+(m+1)\kappa-\delta+{\frac 12} c^{(\mu|\delta,\kappa)})\eta_r\big)$.}


By using any gamma function $G(u|\delta)$ for (a constant multiple of) $\sgm{u}$,
we def\/ine a function $\Phi_{BC}(x;y|\delta,\kappa)$ of Cauchy type,
either by
\begin{gather}\label{eq:phifunBC}
\Phi_{BC}(x;y|\delta,\kappa)
=\dprod{j=1}{m}\dprod{l=1}{n} \dfrac{G(x_j\pm y_l+{\frac 12}(\delta-\kappa)|\delta)}
{G(x_j\pm y_l+{\frac 12}(\delta+\kappa)|\delta)},
\end{gather}
or by
\begin{gather}\label{eq:phifunBC2}
\Phi_{BC}(x;y|\delta,\kappa)
=\dprod{j=1}{m}\dprod{l=1}{n}\,\dprod{\epsilon_1,\epsilon_2=\pm}{}
G(\epsilon_1 x_j+\epsilon_2 y_l+{\tfrac 12}(\delta-\kappa)|\delta).
\end{gather}
Note that the function \eqref{eq:phifunBC2} dif\/fers from \eqref{eq:phifunBC},
only by a multiplicative factor which is $\delta$-periodic in all the
variables $x_j$ ($j=1,\ldots,m$) and $y_l$ ($l=1,\ldots,n$).
We also introduce the function
\begin{gather}\label{eq:psifunBC}
\Psi_{BC}(x;y)=\dprod{j=1}{m}\dprod{l=1}{n}
 \sgm{x_j\pm y_l}
\end{gather}
of dual Cauchy type.

\begin{theorem}\label{thm:BCE}
Suppose that $\sgm{u}$ is a constant multiple of any function in
\eqref{eq:sgmBC}.
\begin{enumerate}\itemsep=0pt
\item[$(1)$] Under the balancing condition
\begin{gather*}
2(m-n)\kappa+c^{(\mu|\delta,\kappa)}=0,
\end{gather*}
the function $\Phi_{BC}(x;y|\delta,\kappa)$
defined as \eqref{eq:phifunBC} or \eqref{eq:phifunBC2}
satisfies the functional equation
\begin{gather*}
E_x^{(\mu|\delta,\kappa)}\Phi_{BC}(x;y|\delta,\kappa)=
E_y^{(\nu|\delta,\kappa)}\Phi_{BC}(x;y|\delta,\kappa),
\end{gather*}
where the parameters $\nu=(\nu_1,\ldots,\nu_{2\rho})$ for the
$y$ variables are defined by
\begin{gather*}
\nu_r={\tfrac 12}(\delta+\kappa)-\mu_r\qquad(r=1,\ldots,2\rho).
\end{gather*}
\item[$(2)$]
Under the balancing condition
\begin{gather*}
2m\kappa+2n\delta+c^{(\mu|\delta,\kappa)}=0,
\end{gather*}
the function $\Psi_{BC}(x;y)$ defined as \eqref{eq:psifunBC} satisfies
the functional equation
\begin{gather*}
\sgm{\kappa} E_x^{(\mu|\delta,\kappa)}\Psi_{BC}(x;y)+
\sgm{\delta} E_y^{(\mu|\kappa,\delta)}\Psi_{BC}(x;y)=0.
\end{gather*}
\end{enumerate}
\end{theorem}
Statement (1) for the cases $m=n$ is due to Ruijsenaars \cite{Ru2005, Ru2009}.
An explicit comparison will be made in Appendix~\ref{appendixB.2} between our
$\Phi_{BC}(x;y|\delta,\kappa)$ and Ruijsenaars' kernel function
of type~$BC$ in~\cite{Ru2009}.

In the case of $0$ variables, the Ruijsenaars operator
$E_{x}^{(\mu|\delta,\kappa)}$ with $m=0$
reduces to the multiplication operator by the constant
\begin{gather*}
C^{(\mu|\delta,\kappa)}
=\dsum{r=1}{\rho}
\dfrac{e(-{\frac 12} c^{(\mu|\delta,\kappa)}\eta_r)
\prod\limits_{s=1}^{2\rho} \sgm{{\frac 12}(\omega_r - \delta) + \mu_s}}
{\sgm{\kappa}
\prod\limits_{s\ne r}
\sgm{{\frac 12}(\omega_r - \omega_s)}
\prod\limits_{s=1}^{\rho}
\sgm{{\frac 12}(\omega_r - \omega_s + \kappa - \delta)}}.
\end{gather*}
We remark that Theorem \ref{thm:BCE}, (2) for $n=0$ implies
that, when $2m\kappa+c^{(\mu|\delta,\kappa)}=0$,
the constant function~$1$ is an eigenfunction for $E_{x}^{(\mu|\delta,\kappa)}$:
\begin{gather*}
\sgm{\kappa} E_{x}^{(\mu|\delta,\kappa)}(1)=
-\sgm{\delta}  C^{(\mu|\kappa,\delta)}.
\end{gather*}
(Statement (1) for $n=0$ gives the same formula since
$C^{(\nu|\delta,\kappa)}=-\sgm{\delta}
C^{(\mu|\kappa,\delta)}/\sgm{\kappa}$.)

In the trigonometric and rational cases, these functions
$\Phi_{BC}(x;y|\delta,\kappa)$ and $\Psi_{BC}(x;y)$ satisfy the following
functional equations without balancing conditions.
\begin{theorem}\label{thm:BCT}
Suppose that $\sgm{u}$ is a constant multiple of
$\sin(\pi u/\omega_1)$ $(\rho=2)$ or
$u$ $(\rho=1)$.
\begin{enumerate}\itemsep=0pt
\item[$(1)$]
The function $\Phi_{BC}(x;y|\delta,\kappa)$
satisfies the functional equation
\begin{gather*}
\sgm{\kappa} E_x^{(\mu|\delta,\kappa)}\Phi_{BC}(x;y|\delta,\kappa)-
\sgm{\kappa} E_y^{(\nu|\delta,\kappa)}\Phi_{BC}(x;y|\delta,\kappa) \\
\qquad{}=
\sgm{2(m - n)\kappa+c^{(\mu|\delta,\kappa)}}
\Phi_{BC}(x;y|\delta,\kappa),
\end{gather*}
where $\nu_r={\frac 12}(\delta+\kappa)-\mu_r$ $(r=1,\ldots,2\rho)$.

\item[$(2)$]
The function $\Psi_{BC}(x;y)$ satisfies
the functional equation
\begin{gather*}
\sgm{\kappa} E_x^{(\mu|\delta,\kappa)}\Psi_{BC}(x;y)+
\sgm{\delta} E_y^{(\mu|\kappa,\delta)}\Psi_{BC}(x;y)
=\sgm{2m\kappa + 2n\delta + c^{(\mu|\delta,\kappa)}} \Psi_{BC}(x;y).
\end{gather*}
\end{enumerate}
\end{theorem}
As a special case $n=0$ of this theorem,
we see that the constant function $1$ is a eigenfunction
of $E_x^{(\mu|\delta,\kappa)}$
for arbitrary values of the parameters.
Theorems \ref{thm:BCE} and \ref{thm:BCT} will
be proved in Section~\ref{section3}.

In the trigonometric and rational $BC$ cases, it is convenient to
introduce another dif\/ference operator
\begin{gather*}
D_{x}^{(\mu|\delta,\kappa)}
=
\dsum{i=1}{m}
A_{i}^{+}(x;\mu|\delta,\kappa)\big(T_{x_i}^{\delta}-1\big)
+
\dsum{i=1}{m}
A_{i}^{-}(x;\mu|\delta,\kappa)\big(T_{x_i}^{-\delta}-1\big)
\end{gather*}
with the same coef\/f\/icients $A_i^{\epsilon}(x;\mu|\delta,\kappa)$
($i=1,\ldots,m; \epsilon=\pm$) as those of the Ruijsenaars opera\-tor~$E_x^{(\mu|\delta,\kappa)}$.
In the trigonometric case, this operator $D_{x}^{(\mu|\delta,\kappa)}$
is a constant multiple of
Koornwinder's $q$-dif\/ference operator expressed in terms
of additive variables.
By using the relation
\begin{gather*}
D_{x}^{(\mu|\delta,\kappa)}=E_x^{(\mu|\delta,\kappa)}-
E_x^{(\mu|\delta,\kappa)}(1),
\end{gather*}
one can rewrite the functional equations in
Theorem \ref{thm:BCT} into those for $D_{x}^{(\mu|\delta,\kappa)}$.

\begin{theorem}\label{thm:BCTD}
Suppose that $\sgm{u}=2\sqrt{-1}\sin(\pi u/\omega_1)$, $\rho=2$
$($resp.~$\sgm{u}=u$, $\rho=1$$)$.
\begin{enumerate}\itemsep=0pt
\item[$(1)$]  The function $\Phi_{BC}(x;y|\delta,\kappa)$
satisfies the functional equation
\begin{gather*}
\sgm{\kappa} D_x^{(\mu|\delta,\kappa)}\Phi_{BC}(x;y|\delta,\kappa)-
\sgm{\kappa} D_y^{(\nu|\delta,\kappa)}\Phi_{BC}(x;y|\delta,\kappa)\\
\qquad{}=
\sgm{m\kappa}\sgm{-n\kappa}\sgm{(m - n)\kappa + c^{(\mu|\delta,\kappa)}}
 \Phi_{BC}(x;y|\delta,\kappa)
\quad(\mbox{resp.} =0),
\end{gather*}
where $\nu_r={\frac 12}(\delta+\kappa)-\mu_r$ $(r=1,\ldots,2\rho)$.

\item[$(2)$] For arbitrary $m$, $n$,
the function $\Psi_{BC}(x;y)$ satisfies
the functional equation
\begin{gather*}
\sgm{\kappa} D_x^{(\mu|\delta,\kappa)}\Psi_{BC}(x;y)+
\sgm{\delta} D_y^{(\mu|\kappa,\delta)}\Psi_{BC}(x;y)
\\
\qquad{}=\sgm{m\kappa}\sgm{n\delta}
\sgm{m\kappa+ n\delta + c^{(\mu|\delta,\kappa)}} \Psi_{BC}(x;y)
\quad(\mbox{resp.} =0).
\end{gather*}
\end{enumerate}
\end{theorem}
Statement (2) of Theorem~\ref{thm:BCTD} recovers
a main result of Mimachi \cite{Mi2001}.
A~proof of Theorem~\ref{thm:BCTD} will be given in Section~\ref{section3.4}.
Also, we will explain in Section~\ref{section4}
how Theorem~\ref{thm:BCTD} is related to the theory of Koornwinder
polynomials.

\section{Derivation of kernel functions}\label{section3}

\subsection{Key identities}\label{section3.1}

We start with a functional identity which decomposes
a product of functions expressed by $\sgm{u}$
into partial fractions.

\begin{proposition}\label{prop:partfrac}
Let $\sgm{u}$ be any nonzero entire function satisfying the Riemann
relation.  Then,
for variables $z$, $(x_1,\ldots, x_N)$, and parameters $(c_1,\ldots,c_N)$,
we have
\begin{gather}\label{eq:partfrac}
\sgm{c} \dprod{j=1}{N}\dfrac{\sgm{z-x_j+c_j}}{\sgm{z-x_j}}
=
\dsum{i=1}{N} \sgm{c_i} \dfrac{\sgm{z-x_i+c}}{\sgm{z-x_i}}
\dprod{1\le j\le N;\,j\ne i}{}\dfrac{\sgm{x_i-x_j+c_j}}{\sgm{x_i-x_j}},
\end{gather}
where $c=c_1+\cdots+c_N$.
\end{proposition}
This identity can be proved by the induction on the number of
factors by using the Riemann relation for $\sgm{u}$.
Note that, by setting $y_j=x_j-c_j$,
this identity can also be rewritten as
\begin{gather*}
\sgm{c}
\dprod{j=1}{N}\dfrac{\sgm{z-y_j}}{\sgm{z-x_j}}
=
\dsum{i=1}{N} \sgm{x_i-y_i} \dfrac{\sgm{z-x_i+c}}{\sgm{z-x_i}}
\dprod{1\le j\le N;\,j\ne i}{}\dfrac{\sgm{x_i-y_j}}{\sgm{x_i-x_j}},
\end{gather*}
where $c=\sum\limits_{j=1}^{N}(x_j-y_j)$.

From this proposition, we obtain the following lemma,
which provides key identities for our proof of kernel relations.

\begin{lemma}\label{lem:keylemma}
Consider $N$ variables $x_1,\ldots,x_N$ and $N$ complex parameters
$c_1,\ldots,c_N$.
\begin{enumerate}\itemsep=0pt
\item[$(1)$]
Let $\sgm{u}$ be any nonzero function satisfying the Riemann relation.
Then under the balancing condition
$\sum\limits_{j=1}^{N}c_j=0$, we have an identity
\begin{gather}\label{eq:eqE}
\dsum{i=1}{N}   \sgm{c_i}
\dprod{1\le j\le N;\,j\ne i}{}\dfrac{\sgm{x_i-x_j+c_j}}{\sgm{x_i-x_j}}
=0.
\end{gather}
as a meromorphic function in $(x_1,\ldots,x_N)$.
\item[$(2)$] Suppose that $\sgm{u}$ is a constant multiple of
$\sin(\pi u/\omega_1)$ or $u$.  Then
for any $c_1,\ldots,c_N\in\mathbb{C}$, we have an identity
\begin{gather}\label{eq:eqTR}
\sum\limits_{i=1}^{N}   \sgm{c_i}
\dprod{1\le j\le N;\,j\ne i}{}\dfrac{\sgm{x_i-x_j+c_j}}{\sgm{x_i-x_j}}
=\left[\sum\limits_{i=1}^{N}c_i\right]
\end{gather}
as a meromorphic function in $(x_1,\ldots,x_N)$.
\end{enumerate}
\end{lemma}

When $\sgm{z}=\sin(\pi z/\omega_1)$, we have
$\sgm{z+a}/\sgm{z+b} \to e(-(a-b)/2\omega_1)$ as
$\mbox{Im}(z/\omega_1)\to \infty$.
This implies in \eqref{eq:partfrac}
\begin{gather*}
\dprod{j=1}{N}
\dfrac{\sgm{z-x_j+c_j}}{\sgm{z-x_j}}\to
e(-c/2\omega_1),\qquad
\dfrac{\sgm{z-x_i+c}}{\sgm{z-x_i}}\to e(-c/2\omega_1)\qquad(i=1,\ldots,m),
\end{gather*}
as $\mbox{Im}(z/\omega_1)\to \infty$.
Hence we obtain \eqref{eq:eqTR} from \eqref{eq:partfrac}.
In the rational case $\sgm{z}=z$, formula
\eqref{eq:eqTR} is derived by a simple limiting procedure
$z\to\infty$.

\subsection[Case of type $A$]{Case of type $\boldsymbol{A}$}\label{section3.2}

We apply the key identities~\eqref{eq:eqE} and~\eqref{eq:eqTR}
for studying kernel functions.
For two sets of variables
$x=(x_1,\ldots,x_m)$ and $y=(y_1,\ldots,y_n)$,
we consider the following meromorphic function:
\begin{gather*}
F(z)=
\dprod{j=1}{m}\dfrac{\sgm{z-x_j+\kappa}}{\sgm{z-x_j}}
\dprod{l=1}{n}\dfrac{\sgm{z-y_l+v+\lambda}}{\sgm{z-y_l+v}}
\end{gather*}
where $\kappa$, $\lambda$ and $v$ are complex parameters.
Then by Proposition \ref{prop:partfrac} this function $F(z)$ is expanded
as
\begin{gather*}
\sgm{c} F(z) =
\sgm{\kappa} \dsum{i=1}{m}
\dfrac{\sgm{z-x_i+c}}{\sgm{z-x_i}}
\dprod{1\le j\le m;\,j\ne i}{}
\dfrac{\sgm{x_i-x_j+\kappa}}{\sgm{x_i-x_j}}
\dprod{l=1}{n}
\dfrac{\sgm{x_i-y_l+v+\lambda}}{\sgm{x_i-y_l+v}}
\nonumber\\
\phantom{\sgm{c} F(z) =}{} +\sgm{\lambda}
\dsum{k=1}{n}
\dfrac{\sgm{z-y_k+v+c}}{\sgm{z-y_k+v}}
\dprod{1\le l\le n;\,l\ne k}{}
\dfrac{\sgm{y_k-y_l+\lambda}}{\sgm{y_k-y_l}}
\dprod{j=1}{m}
\dfrac{\sgm{y_k-x_j-v+\kappa}}{\sgm{y_k-x_j-v}}
\nonumber\\
\phantom{\sgm{c} F(z)}{}=\sgm{\kappa} \dsum{i=1}{m}
\dfrac{\sgm{z-x_i+c}}{\sgm{z-x_i}}
A_i(x;\kappa)
\dprod{l=1}{n}
\dfrac{\sgm{x_i-y_l+v+\lambda}}{\sgm{x_i-y_l+v}}
\nonumber\\
\phantom{\sgm{c} F(z)=}{}
+\sgm{\lambda}
\dsum{k=1}{n}
\dfrac{\sgm{z-y_k+v+c}}{\sgm{z-y_k+v}}
A_k(y;\lambda)
\dprod{j=1}{m}
\dfrac{\sgm{y_k-x_j-v+\kappa}}{\sgm{y_k-x_j-v}},
\end{gather*}
where $c=m\kappa+n\lambda$,
with the coef\/f\/icients $A_i(x;\kappa)$, $A_k(y;\lambda)$
of Ruijsenaars operators in $x$ variables and $y$ variables.
By Lemma~\ref{lem:keylemma},
under the balancing condition $c=m\kappa+n\lambda=0$,
we have
\begin{gather*}
\sgm{\kappa} \dsum{i=1}{m}
A_i(x;\kappa)
\dprod{l=1}{n}
\dfrac{\sgm{x_i-y_l+v+\lambda}}{\sgm{x_i-y_l+v}}
+\sgm{\lambda}
\dsum{k=1}{n}
A_k(y;\lambda)
\dprod{j=1}{m}
\dfrac{\sgm{y_k-x_j-v+\kappa}}{\sgm{y_k-x_j-v}}
=0.
\end{gather*}
Also,
when $\sgm{u}$ is a constant multiple of $\sin(\pi u/\omega_1)$
or $u$,
we always have
\begin{gather*}
\sgm{\kappa} \dsum{i=1}{m}
A_i(x;\kappa)
\dprod{l=1}{n}
\dfrac{\sgm{x_i-y_l+v+\lambda}}{\sgm{x_i-y_l+v}}
+\sgm{\lambda}
\dsum{k=1}{n}
A_k(y;\lambda)
\dprod{j=1}{m}
\dfrac{\sgm{y_k-x_j-v+\kappa}}{\sgm{y_k-x_j-v}}
=\sgm{m\kappa+n\lambda}.
\end{gather*}

We now try to f\/ind a function $\Phi(x;y)$
satisfying the system
of f\/irst-order dif\/ference equations
\begin{gather}
T_{x_i}^{\delta}\Phi(x;y)
=
\dprod{l=1}{n}
\dfrac{\sgm{x_i-y_l+v+\lambda}}{\sgm{x_i-y_l+v}}
\Phi(x;y)\qquad(i=1,\ldots,m),\nonumber\\
T_{y_k}^{\tau}\Phi(x;y)
=
\dprod{j=1}{m}
\dfrac{\sgm{y_k-x_j-v+\kappa}}{\sgm{y_k-x_j-v}}
\Phi(x;y)\qquad (k=1,\ldots,n),\label{eq:LPA}
\end{gather}
where $\tau$ stands for the unit scale of dif\/ference operators for $y$
variables.
In order to f\/ix the idea, assume that the balancing condition
$m\kappa+n\lambda=0$ is satisf\/ied.
Then,
any solution of this system should satisfy the functional equation
\begin{gather*}
\sgm{\kappa} \dsum{i=1}{m}
A_i(x;\kappa) T_{x_i}^{\delta}\Phi(x;y)
+\sgm{\lambda}
\dsum{k=1}{n}
A_k(y;\lambda)
T_{y_k}^{\tau}\Phi(x;y)
=0,
\end{gather*}
namely,
\begin{gather*}
\sgm{\kappa} D_{x}^{(\delta,\kappa)}\Phi(x;y)+
\sgm{\lambda} D_{y}^{(\tau,\lambda)}\Phi(x;y)=0.
\end{gather*}

The compatibility condition for the system \eqref{eq:LPA} of
dif\/ference equations is given by
\begin{gather*}
\dfrac{\sgm{x_i - y_k - \tau + v + \lambda}\sgm{x_i - y_k + v}}
{\sgm{x_i - y_k - \tau + v}\sgm{x_i - y_k + v + \lambda}}
=
\dfrac{\sgm{y_k - x_i - \delta - v + \kappa}\sgm{y_k - x_i - v}}
{\sgm{y_k - x_i - \delta - v}\sgm{y_k - x_i - v + \kappa}}
\end{gather*}
for $i=1,\ldots,m$ and $k=1,\ldots,n$.
From this we see there are (at least) two cases where the dif\/ference
equation \eqref{eq:LPA} becomes compatible:
\begin{alignat*}{5}
& \mbox{(case 1) :}\qquad&&
\tau=-\delta,\qquad && \lambda=-\kappa,\qquad && v:\ \mbox{arbitrary},&
\\
& \mbox{(case 2)  :}\qquad&&
\tau=\kappa,\qquad && \lambda=\delta,\qquad && v: \ \mbox{arbitrary}.&
\end{alignat*}
In the f\/irst case, the dif\/ference equation to be solved is:
\begin{gather*}
T_{x_i}^{\delta}\Phi(x;y)
 =
\dprod{l=1}{n}
\dfrac{\sgm{x_i-y_l+v-\kappa}}{\sgm{x_i-y_l+v}}
\Phi(x;y)\qquad (i=1,\ldots,m),\\
T_{y_k}^{-\delta}\Phi(x;y)
 =
\dprod{j=1}{m}
\dfrac{\sgm{y_k-x_j-v+\kappa}}{\sgm{y_k-x_j-v}}
\Phi(x;y)\qquad (k=1,\ldots,n).
\end{gather*}
This system is solved by
\begin{gather*}
\Phi_{A}^{-}(x;y|\delta,\kappa)=
\dprod{j=1}{m}\dprod{l=1}{n}
\dfrac{G(x_j-y_l+v-\kappa|\delta)}{G(x_j-y_l+v|\delta)}.
\end{gather*}
Hence we see that $\Phi_{A}^{-}(x;y|\delta,\kappa)$
satisf\/ies the functional equation
\begin{gather*}
D_{x}^{(\delta,\kappa)}\Phi_{A}^{-}(x;y|\delta,\kappa)
-D_{y}^{(-\delta,-\kappa)}\Phi_{A}^{-}(x;y|\delta,\kappa)=0,
\end{gather*}
under the balancing condition $(m-n)\kappa=0$, namely, $m=n$.
By setting $\Phi_{A}(x;y)=\Phi_{A}^{-}(x;-y)$, we
obtain
\begin{gather*}
D_{x}^{(\delta,\kappa)}\Phi_{A}(x;y|\delta,\kappa)
-D_{y}^{(\delta,\kappa)}\Phi_{A}(x;y|\delta,\kappa)=0
\end{gather*}
as in Theorem \ref{thm:AE}.
The second case
\begin{gather*}
\arraycolsep=2pt
\begin{array}{rll}
T_{x_i}^{\delta}\Phi(x;y)
&=
\dprod{l=1}{n}
\dfrac{\sgm{x_i-y_l+v+\delta}}{\sgm{x_i-y_l+v}}
\Phi(x;y)\quad&(i=1,\ldots,m),\\[14pt]
T_{y_k}^{\kappa}\Phi(x;y)
&=
\dprod{j=1}{m}
\dfrac{\sgm{y_k-x_j-v+\kappa}}{\sgm{y_k-x_j-v}}
\Phi(x;y)\quad&(k=1,\ldots,n)
\end{array}
\end{gather*}
is solved by
\begin{gather*}
\Psi_{A}(x;y)=
\dprod{j=1}{m}\dprod{l=1}{n} \sgm{x_j-y_l+v}.
\end{gather*}
Hence, under the balancing condition
$m\kappa+n\delta=0$, we have
\begin{gather*}
\sgm{\kappa} D_{x}^{(\delta,\kappa)}\Psi_{A}(x;y)+
\sgm{\delta} D_{y}^{(\kappa,\delta)}\Psi_{A}(x;y)=0.
\end{gather*}
This completes the proof of Theorem~\ref{thm:AE}.

When $\sgm{u}$ is a constant multiple of $\sin(\pi u/\omega_1)$
or $u$, for any solution $\Phi(x;y)$
of the sys\-tem~\eqref{eq:LPA} of dif\/ference equations
we have
\begin{gather*}
\sgm{\kappa} \dsum{i=1}{m}
A_i(x;\kappa) T_{x_i}^{\delta}\Phi(x;y)
+\sgm{\lambda}
\dsum{k=1}{n}
A_k(y;\lambda)
T_{y_k}^{\kappa}\Phi(x;y)
=\sgm{m\kappa+n\lambda}\Phi(x;y),
\end{gather*}
namely,
\begin{gather*}
\sgm{\kappa} D_{x}^{(\delta,\kappa)}\Phi(x;y)+
\sgm{\lambda} D_{y}^{(\tau,\lambda)}\Phi(x;y)=\sgm{m\kappa+n\lambda}\Phi(x;y)
\end{gather*}
without imposing the balancing condition.
This implies that the functions $\Phi_{ A}(x;y|\delta,\kappa)$
and $\Psi_{ A}(x;y)$
in these trigonometric and rational cases
satisfy the functional equations
\begin{gather*}
D_{x}^{(\delta,\kappa)}\Phi_{ A}(x;y|\delta,\kappa)
- D_{y}^{(\delta,\kappa)}\Phi_{ A}(x;y|\delta,\kappa)=
\dfrac{\sgm{(m-n)\kappa}}{\sgm{\kappa}}
\Phi_{ A}(x;y|\delta,\kappa)
\end{gather*}
and
\begin{gather*}
\sgm{\kappa} D_{x}^{(\delta,\kappa)}\Psi_{ A}(x;y)+
\sgm{\delta} D_{y}^{(\kappa,\delta)}\Psi_{ A}(x;y)
=\sgm{m\kappa+n\delta} \Psi_{ A}(x;y).
\end{gather*}
respectively, as stated in Theorem \ref{thm:AT}.

\subsection[Case of type $BC$]{Case of type $\boldsymbol{BC}$}\label{section3.3}

In this $BC$ case, we assume that $\sgm{u}$ is a constant multiple
of one of the following functions:
\begin{gather*}
u\ \ (\rho=1),\qquad
\sin(\pi u/\omega_1)\ \ (\rho=2),\qquad
e\left(au^2\right)\sigma(u;\Omega)\ \ (\rho=4).
\end{gather*}
In order to discuss dif\/ference operators of type $BC$,
we consider the meromorphic function
\begin{gather}\label{eq:defF}
F(z)=
\dfrac{\prod\limits_{s=1}^{2\rho}\sgm{z+\mu_s}}
{\prod\limits_{s=1}^{\rho}\sgm{z + {\frac 12}(\delta - \omega_s)}
\sgm{z + {\frac 12}(\kappa - \omega_s)}}
\dprod{j=1}{m}
\dfrac{\sgm{z\pm x_j+\kappa}}{\sgm{z\pm x_j}}
\dprod{l=1}{n}
\dfrac{\sgm{z\pm y_l+v+\lambda}}{\sgm{z\pm y_l+v}}.
\end{gather}
By Proposition~\ref{prop:partfrac} it can be expanded as
\begin{gather}
\sgm{c} F(z)=
\dsum{i=1}{m}\dfrac{\sgm{z-x_i+c}}{\sgm{z-x_i}}  P_i^{+}
+
\dsum{i=1}{m}\dfrac{\sgm{z+x_i+c}}{\sgm{z+x_i}}  P_i^{-}
\nonumber\\
\phantom{\sgm{c} F(z)=}{}
+\dsum{k=1}{n}\dfrac{\sgm{z-y_k+v+c}}{\sgm{z-y_k+v}}  Q_k^{+}
+
\dsum{k=1}{n}\dfrac{\sgm{z+y_k+v+c}}{\sgm{z+y_k+v}}  Q_k^{-}\nonumber\\
\phantom{\sgm{c} F(z)=}{}
 +
\dsum{r=1}{\rho}\dfrac{\sgm{z+{\frac 12}(\delta-\omega_r)+c}}
{\sgm{z+{\frac 12}(\delta-\omega_r)}}  R_r
+
\dsum{r=1}{\rho}\dfrac{\sgm{z+{\frac 12}(\kappa-\omega_r)+c}}
{\sgm{z+{\frac 12}(\kappa-\omega_r)}}  S_r\label{eq:expF}
\end{gather}
into partial fractions, where
\begin{gather*}
c=2m\kappa+2n\lambda+\sum\limits_{s=1}^{2\rho}\mu_s
-\tfrac{\rho}{2}(\delta+\kappa)
+\sum\limits_{s=1}^{\rho}\omega_s=2m\kappa+2n\lambda+c^{(\mu|\delta,\kappa)}.
\end{gather*}
Also by Lemma \ref{lem:keylemma}, we see that
expression
\begin{gather*}
\dsum{1\le i\le m;\,\epsilon=\pm}{}P_i^\epsilon+
\dsum{1\le k\le n;\,\epsilon=\pm}{}Q_k^\epsilon+
\dsum{1\le r\le \rho}{}R_r+
\dsum{1\le r\le \rho}{}S_r
\end{gather*}
reduces to 0 when the balancing condition
$c=0$ is satisf\/ied, or to
$\sgm{c}$ when $\sgm{u}$ is trigonometric or rational.
A remarkable fact is that,
if the parameter $v$ is chosen appropriately,
then the expansion coef\/f\/icients of $F(z)$ are expressed in terms
of the coef\/f\/icients of Ruijsenaars operators of type~$BC$.

\begin{proposition}\label{prop:expF}
When $v={\frac 12}(\delta-\lambda)$, the function $F(z)$
defined by \eqref{eq:defF} is expressed as follows
in terms of the coefficients of Ruijsenaars operators:
\begin{gather*}
\sgm{c} F(z)
=
\dfrac{\sgm{c} \prod\limits_{s=1}^{2\rho}\sgm{z+\mu_s}}
{\prod\limits_{s=1}^{\rho}\sgm{z + {\frac 12}(\delta - \omega_s)}
\sgm{z + {\frac 12}(\kappa - \omega_s)}}
 \dprod{j=1}{m}
\dfrac{\sgm{z\pm x_j+\kappa}}{\sgm{z\pm x_j}}
 \dprod{l=1}{n}
\dfrac{\sgm{z\pm y_l+v+\lambda}}{\sgm{z\pm y_l+v}}
\\
\phantom{\sgm{c} F(z) }{} =
\sgm{\kappa}
 \dsum{1\le i\le m;\,\epsilon=\pm}{}
\dfrac{\sgm{z-\epsilon x_i+c}}
{\sgm{z-\epsilon x_i}}  A_i^{\epsilon}(x;\mu|\delta,\kappa)
 \dprod{l=1}{n}
\dfrac{\sgm{\epsilon x_i\pm y_l+{\frac 12}(\delta + \lambda)}}
{\sgm{\epsilon x_i\pm y_l+{\frac 12}(\delta - \lambda)}}
\\
\phantom{\sgm{c} F(z) =}{} +\sgm{\kappa}
 \dsum{r=1}{\rho}
e\big({\tfrac 12} c \eta_r\big)
\dfrac{\sgm{z+{\frac 12}(\delta-\omega_r)+c}}
{\sgm{z+{\frac 12}(\delta-\omega_r)}}
A_r^{0}(x;\mu|\delta,\kappa)
\\
\phantom{\sgm{c} F(z) =}{}+\sgm{\lambda}
\dsum{1\le k\le n;\, \epsilon=\pm}{}
\dfrac{\sgm{z-\epsilon y_k+v+c}}{\sgm{z-\epsilon y_k+v}}
A_k^{\epsilon}(y;\nu|\tau,\lambda)
\dprod{j=1}{m}
\dfrac{\sgm{\epsilon y_k\pm x_j+{\frac 12}(\tau+\kappa)}}
{\sgm{\epsilon y_k\pm x_j+{\frac 12}(\tau-\kappa)}}
\\
\phantom{\sgm{c} F(z) =}{}+
\sgm{\lambda}
\dsum{r=1}{\rho}
e\big({\tfrac 12} c \eta_r\big)
\dfrac{\sgm{z+{\frac 12}(\kappa-\omega_r)+c}}
{\sgm{z+{\frac 12}(\kappa-\omega_r)}} A_r^{0}(y;\nu|\tau,\lambda),
\end{gather*}
where
\begin{gather*}
c=
2m\kappa+2n\lambda+c^{(\mu|\delta,\kappa)},
\qquad
\tau=\kappa+\lambda-\delta, \qquad
\nu=(\mu_1-v,\ldots,\mu_{2\rho}-v).
\end{gather*}
\end{proposition}

\begin{proof}
The expansion coef\/f\/icients in \eqref{eq:expF}
are determined from the residues at
the corresponding poles.
We f\/irst remark that
\begin{gather*}
P_i^{+}=
\dfrac{\sgm{\kappa} \prod\limits_{s=1}^{2\rho}\sgm{x_i+\mu_s}}
{\prod\limits_{s=1}^{\rho}\sgm{x_i + {\frac 12}(\delta - \omega_s)}
\sgm{x_i + {\frac 12}(\kappa - \omega_s)}}
\dfrac{\sgm{2x_i+\kappa}}{\sgm{2x_i}}
\dprod{j\ne i}{m}
\dfrac{\sgm{x_i\pm x_j+\kappa}}{\sgm{x_i\pm x_j}}
\dprod{l=1}{n}
\dfrac{\sgm{x_i\pm y_l+v+\lambda}}{\sgm{x_i\pm y_l+v}}.
\end{gather*}
Since
\begin{gather*}
\dfrac{\sgm{2x_i+\kappa}}{\sgm{2x_i}}
=\dprod{s=1}{\rho}
\dfrac{\sgm{x_i+{\frac 12}(\kappa-\omega_s)}}{\sgm{x_i-{\frac 12}\omega_s}},
\end{gather*}
we obtain
\begin{gather*}
P_i^{+} =
\dfrac{\sgm{\kappa} \prod\limits_{s=1}^{2\rho}\sgm{x_i+\mu_s}}
{\prod\limits_{s=1}^{\rho}\sgm{x_i-{\frac 12}\omega_s} \sgm{x_i+{\frac 12}(\delta-\omega_s)}}
\dprod{j\ne i}{m}
\dfrac{\sgm{x_i\pm x_j+\kappa}}{\sgm{x_i\pm x_j}}
\dprod{l=1}{n}
\dfrac{\sgm{x_i\pm y_l+v+\lambda}}{\sgm{x_i\pm y_l+v}}
\\
 \phantom{P_i^{+}}{} =\sgm{\kappa} A_i^{+}(x;\mu|\delta,\kappa)
\dprod{l=1}{n}
\dfrac{\sgm{x_i\pm y_l+v+\lambda}}{\sgm{x_i\pm y_l+v}}.
\end{gather*}
Note that $P_i^{-}$ is obtained from $P_i^{+}$ by replacing $x_j$ with
$-x_j$ ($j=1,\ldots,m$).
We next look at the coef\/f\/icient $Q_k^{+}$:
\begin{gather*}
Q_k^{+} =
\dfrac{\sgm{\lambda} \prod\limits_{s=1}^{2\rho}\sgm{y_k+\mu_s-v}}
{\prod\limits_{s=1}^{\rho}\sgm{y_k + {\frac 12}(\delta - \omega_s)-v}
\sgm{y_k + {\frac 12}(\kappa - \omega_s)-v}}
\dfrac{\sgm{2y_k+\lambda}}{\sgm{2y_k}}
\\
\phantom{Q_k^{+} =}{}\times
\dprod{l\ne k}{}
\dfrac{\sgm{y_k\pm y_l+\lambda}}{\sgm{y_k\pm y_l}}
\dprod{j=1}{m}
\dfrac{\sgm{y_k\pm x_j+\kappa-v}}{\sgm{y_k\pm x_j-v}}.
\end{gather*}
In view of
\begin{gather*}
\dfrac{\sgm{2y_k+\lambda}}{\sgm{2y_k}}
=\dprod{s=1}{\rho}
\dfrac{\sgm{y_k+{\frac 12}(\lambda-\omega_s)}}{\sgm{y_k-{\frac 12}\omega_s}}
\end{gather*}
we set ${\frac 12}\delta-v={\frac 12}\lambda$, namely, $v={\frac 12}(\delta-\lambda)$.
Then we have
\begin{gather*}
Q_k^{+} =
\dfrac{\sgm{\lambda} \prod\limits_{s=1}^{2\rho}\sgm{y_k + \mu_s - v}}
{\prod\limits_{s=1}^{\rho}\sgm{y_k - {\frac 12}\omega_s}
\sgm{y_k + {\frac 12}(\tau - \omega_s)}}
\dprod{l\ne k}{}
\dfrac{\sgm{y_k\pm y_l + \lambda}}{\sgm{y_k\pm y_l}}
\dprod{j=1}{m}
\dfrac{\sgm{y_k\pm x_j + \kappa - v}}{\sgm{y_k\pm x_j - v}}\\
\phantom{Q_k^{+}}{} =
\sgm{\lambda}
A_{k}^{+}(y;\nu|\tau,\lambda)
\dprod{j=1}{m}
\dfrac{\sgm{y_k\pm x_j + \kappa - v}}{\sgm{y_k\pm x_j - v}},\\
  v={\tfrac 12}(\delta-\lambda),\qquad \tau=\kappa+\lambda-\delta,\qquad
\nu=(\mu_1-v,\ldots,\mu_{2\rho}-v).
\end{gather*}
The coef\/f\/icient $Q_k^{-}$ is obtained from $Q_k^{+}$ by the sign change
$y_l\to -y_l$.
The coef\/f\/icient $R_r$ is given~by
\begin{gather*}
R_r =
\dfrac{
\prod\limits_{s=1}^{2\rho}\sgm{{\frac 12}(\omega_r - \delta)+\mu_s}}
{\prod\limits_{s\ne r}\sgm{{\frac 12}\omega_{rs}}
 \prod\limits_{s=1}^{\rho}
\sgm{{\frac 12}\omega_{rs} + {\frac 12}(\kappa-\delta)}}
\dprod{j=1}{m}
\dfrac{\sgm{{\frac 12}(\omega_r - \delta)\pm x_j + \kappa}}
{\sgm{{\frac 12}(\omega_r - \delta)\pm x_j}}
\dprod{l=1}{n}
\dfrac{\sgm{{\frac 12}(\omega_r - \delta)\pm y_l +v+ \lambda}}
{\sgm{{\frac 12}(\omega_r - \delta)\pm y_l+v}},
\end{gather*}
where $\omega_{rs}=\omega_r-\omega_s$.
When $v={\frac 12}(\delta-\lambda)$,
we have
\begin{gather*}
\dfrac{\sgm{{\frac 12}(\omega_r - \delta)\pm y_l+v+\lambda}}
{\sgm{{\frac 12}(\omega_r - \delta)\pm y_l+v}}
=
\dfrac{\sgm{{\frac 12}(\omega_r + \lambda)\pm y_l}}
{\sgm{{\frac 12}(\omega_r - \lambda)\pm y_l}}
=e(\lambda\eta_r)
\end{gather*}
by the quasi-periodicity of $\sgm{u}$.
In this way $y$ variables disappear from $R_r$:
\begin{gather*}
R_r=
\dfrac{
e(n\lambda \eta_r)
\prod\limits_{s=1}^{2\rho}\sgm{{\frac 12}(\omega_r - \delta) + \mu_s}}
{\prod\limits_{s\ne r} \sgm{{\frac 12}\omega_{rs}}
\prod\limits_{s=1}^{\rho}
\sgm{{\frac 12}\omega_{rs}+ {\frac 12}(\kappa-\delta)}}
\dprod{j=1}{m}
\dfrac{\sgm{{\frac 12}(\omega_r- \delta)\pm x_j+\kappa}}
{\sgm{{\frac 12}(\omega_r - \delta)\pm x_j}}.
\end{gather*}
Note that the exponential factor $e(n\lambda \eta_r)$ is
nontrivial only in the elliptic case.  In any case,
from $2n\lambda=c-2m\kappa-c^{(\mu|\delta,\kappa)}$ we have
\begin{gather*}
R_r=\sgm{\kappa} e\big({\tfrac 12} c \eta_r\big)  A_r^{0}(x;\mu|\delta,\kappa).
\end{gather*}
Finally, when $v={\frac 12}(\delta-\lambda)$ and $\tau=\kappa+\lambda-\delta$,
we can rewrite $S_r$ as
\begin{gather*}
S_r =
\dfrac{\prod\limits_{s=1}^{2\rho}\sgm{{\frac 12}(\omega_r - \kappa)+\mu_s}}
{\prod\limits_{s\ne r}
\sgm{{\frac 12}\omega_{rs}}
\prod\limits_{s=1}^{\rho}
\sgm{{\frac 12}\omega_{rs} +
{\frac 12}(\delta - \kappa)}}
\dprod{j=1}{m}
\dfrac{\sgm{{\frac 12}(\omega_r + \kappa)\pm x_j}}
{\sgm{{\frac 12}(\omega_r - \kappa)\pm x_j}}
\dprod{l=1}{n}
\dfrac{\sgm{{\frac 12}(\omega_r - \kappa)\pm y_l + v + \lambda}}
{\sgm{{\frac 12}(\omega_r - \kappa)\pm y_l + v}}
\\
\phantom{S_r}{} =
\dfrac{e(m\kappa \eta_r)
\prod\limits_{s=1}^{8}\sgm{{\frac 12}(\omega_r - \tau)+\mu_s-v}}
{\prod\limits_{s\ne r}
\sgm{{\frac 12}\omega_{rs}}
\prod\limits_{s=1}^{4}
\sgm{{\frac 12}\omega_{rs} +
{\frac 12}(\lambda - \tau)}}
\dprod{l=1}{n}
\dfrac{\sgm{{\frac 12}(\omega_r - \tau)\pm y_l + \lambda}}
{\sgm{{\frac 12}(\omega_r - \tau)\pm y_l}}.
\end{gather*}
Note that $x$ variables have disappeared again
by the quasi-periodicity of $\sgm{u}$.
Since $2m\kappa=c-2n\lambda-c^{(\mu|\delta,\kappa)}$ and
$c^{(\mu|\delta,\kappa)}=c^{(\nu|\tau,\lambda)}$,
we have
\begin{gather*}
S_r=\sgm{\lambda}e\big({\tfrac 12} c \eta_r\big)  A_r^{0}(y;\nu|\tau,\lambda).
\end{gather*}
This completes the proof of proposition.
\end{proof}

In what follows we set $v={\frac 12}(\delta-\lambda)$ and
$\tau=\kappa+\lambda-\delta$.
Then Proposition~\ref{prop:expF}, combined with
Lemma~\ref{lem:keylemma}, implies
\begin{gather*}
\sgm{\kappa}
\dsum{1\le i\le m;\ \epsilon=\pm}{}
A_i^{\epsilon}(x;\mu|\delta,\kappa)
\dprod{l=1}{n}
\dfrac{\sgm{\epsilon x_i\pm y_l+{\frac 12}(\delta+\lambda)}}
{\sgm{\epsilon x_i\pm y_l+{\frac 12}(\delta-\lambda)}}
+\sgm{\kappa}
\dsum{s=1}{\rho}
A_r^{0}(x;\mu|\delta,\kappa)
\nonumber\\
\qquad{}
+\sgm{\lambda}
\dsum{1\le k\le n;\ \epsilon=\pm}{}
A_k^{\epsilon}(y;\nu|\tau,\lambda)
 \dprod{j=1}{m}
\dfrac{\sgm{\epsilon y_k\pm x_j+{\frac 12}(\tau + \kappa)}}
{\sgm{\epsilon y_k\pm x_j+{\frac 12}(\tau - \kappa)}}
+
\sgm{\lambda}
\dsum{r=1}{\rho}
A_r^{0}(y;\nu|\tau,\lambda)
\nonumber\\
\quad{} =0,\quad(\mbox{resp.}\
=\sgm{2m\kappa+2n\lambda+c^{(\mu|\delta,\kappa)}}),
\end{gather*}
when the balancing condition
$2m\kappa+2n\lambda+c^{(\mu|\delta,\kappa)}=0$
is satisf\/ied (resp.~when $\sgm{u}$ is trigonometric or rational.)
Hence we have

\begin{proposition}\label{prop:LPBC}
Suppose that the parameters $\delta$, $\kappa$, $\tau$, $\lambda$ satisfy
the relation $\delta+\tau=\kappa+\lambda$,
and define $\nu=(\nu_1,\ldots,\nu_{2\rho})$ by
\begin{gather*}
\nu_s=\mu_s-{\tfrac 12}(\delta-\lambda)=\mu_s+{\tfrac 12}(\tau-\kappa)
\qquad(s=1,\ldots,2\rho).
\end{gather*}
Let $\Phi(x;y)$ any meromorphic function in the variables
$x=(x_1,\ldots,x_m)$ and $y=(y_1,\ldots,y_n)$ satisfying
the system of first-order difference equations
\begin{gather}
T_{x_i}^{\delta}\Phi(x;y)
 =
\dprod{l=1}{n}
\dfrac{\sgm{x_i\pm y_l+{\frac 12}(\delta + \lambda)}}
{\sgm{x_i\pm y_l+{\frac 12}(\delta - \lambda)}}
\Phi(x;y)\qquad (i=1,\ldots,m),
\nonumber\\
T_{y_k}^{\tau}\Phi(x;y)
 =
\dprod{j=1}{m}
\dfrac{\sgm{y_k\pm x_j+{\frac 12}(\tau + \kappa)}}
{\sgm{y_k\pm x_j+{\frac 12}(\tau - \kappa)}}
\Phi(x;y)
\qquad (k=1,\ldots,n).\label{eq:LP}
\end{gather}

\begin{enumerate}\itemsep=0pt
\item[$(1)$] If the balancing condition
$2m\kappa+2n\lambda+c^{(\mu|\delta,\kappa)}=0$ is
satisfied,
then
$\Phi(x;y)$ satisfies the functional equation
\begin{gather*}
\sgm{\kappa} E_{x}^{(\mu|\delta,\kappa)}\Phi(x;y)
+
\sgm{\lambda} E_{y}^{(\nu|\tau,\lambda)}\Phi(x;y)
=0.
\end{gather*}
\item[$(2)$]  If $\sgm{u}$ is a constant multiple
of $\sin(\pi u/\omega_1)$ or $u$, then
$\Phi(x;y)$ satisfies the functional equation
\begin{gather*}
\sgm{\kappa} E_{x}^{(\mu|\delta,\kappa)}\Phi(x;y)
+
\sgm{\lambda} E_{y}^{(\nu|\tau,\lambda)}\Phi(x;y)
=\sgm{2m\kappa+2n\lambda+c^{(\mu|\delta,\kappa)}}\Phi(x;y).
\end{gather*}
\end{enumerate}
\end{proposition}

In fact there are essentially two cases
where the system of f\/irst-order linear dif\/ference equa\-tions~\eqref{eq:LP} become compatible:{\samepage
\begin{alignat}{6}
& \mbox{(case 1):}\qquad && \tau=-\delta,\qquad  && \lambda=-\kappa, \qquad  && \nu_s=\mu_s-{\tfrac 12}(\delta+\kappa) \qquad &&(s=1,\ldots,2\rho);& \nonumber
\\
& \mbox{(case 2):}\qquad && \tau=\kappa,\qquad &&\lambda=\delta, \qquad && \nu_s=\mu_s\qquad && (s=1,\ldots,2\rho).& \label{eq:twocasesBC}
\end{alignat}}

\noindent
In the f\/irst case, the system \eqref{eq:LP}
of dif\/ference equations to be solved is:
\begin{gather*}
T_{x_i}^{\delta}\Phi(x;y)
=
\dprod{l=1}{n}
\dfrac{\sgm{x_i\pm y_l+{\frac 12}(\delta - \kappa)}}
{\sgm{x_i\pm y_l+{\frac 12}(\delta + \kappa)}}
\Phi(x;y)\qquad (i=1,\ldots,m),
\nonumber\\
T_{y_k}^{\delta}\Phi(x;y)
 =
\dprod{j=1}{m}
\dfrac{\sgm{y_k\pm x_j+{\frac 12}(\delta - \kappa)}}
{\sgm{y_k\pm x_j+{\frac 12}(\delta + \kappa)}}
\Phi(x;y)
\qquad (k=1,\ldots,n).
\end{gather*}
It is solved either by the function
\begin{gather*}
\Phi_{BC}(x;y|\delta,\kappa)
=\dprod{j=1}{m}\dprod{l=1}{n}
\dfrac{G(x_j\pm y_l+{\frac 12}(\delta-\kappa)|\delta)}
{G(x_j\pm y_l+{\frac 12}(\delta+\kappa)|\delta)},
\end{gather*}
or by
\begin{gather*}
\Phi_{BC}(x;y|\delta,\kappa)
=\dprod{j=1}{m}\dprod{l=1}{n} \dprod{\epsilon_1,\epsilon_2=\pm}{}
G(\epsilon_1 x_j+\epsilon_2 y_l+{\tfrac 12}(\delta-\kappa)|\delta).
\end{gather*}
Hence we see that,
under the balancing condition
$2(m-n)\kappa+c^{(\mu|\delta,\kappa)}=0$,
$\Phi_{ BC}(x;y|\delta,\kappa)$ satisf\/ies the functional equation
\begin{gather*}
E_{x}^{(\mu|\delta,\kappa)}\Phi_{ BC}(x;y|\delta,\kappa)
-E_{y}^{(\nu|-\delta,-\kappa)}\Phi_{ BC}(x;y|\delta,\kappa)=0,
\end{gather*}
for $\nu_s=\mu_s-{\frac 12}(\delta+\kappa)$
($s=1,\ldots,2\rho$).
By symmetry \eqref{eq:symE} with respect to the sign change
of parameters, we obtain
\begin{gather*}
E_{x}^{(\mu|\delta,\kappa)}\Phi_{ BC}(x;y|\delta,\kappa)
-E_{y}^{(\nu|\delta,\kappa)}\Phi_{ BC}(x;y|\delta,\kappa)=0,
\end{gather*}
for $\nu_s={\frac 12}(\delta+\kappa)-\mu_s$
($s=1,\ldots,2\rho$).
The second case
\begin{gather*}
T_{x_i}^{\delta}\Phi(x;y)
 =
\dprod{l=1}{n}
\dfrac{\sgm{x_i\pm y_l+\delta}}
{\sgm{x_i\pm y_l}}
\Phi(x;y)\qquad (i=1,\ldots,m),
\\
T_{y_k}^{\kappa}\Phi(x;y)
 =
\dprod{j=1}{m}
\dfrac{\sgm{y_k\pm x_j+\kappa}}
{\sgm{y_k\pm x_j}}
\Phi(x;y)
\qquad (k=1,\ldots,n)
\end{gather*}
is solved by
\begin{gather*}
\Psi_{BC}(x;y)=
\dprod{j=1}{m}\dprod{l=1}{n} \sgm{x_j\pm y_l}.
\end{gather*}
Hence, under the balancing condition
$2m\kappa+2n\delta+c^{(\mu|\delta,\kappa)}=0$,
this function satisf\/ies the functional equation
\begin{gather*}
\sgm{\kappa} E_x^{(\mu|\delta,\kappa)} \Psi_{BC}(x;y)+
\sgm{\delta} E_y^{(\mu|\kappa,\delta)} \Psi_{BC}(x;y)=0.
\end{gather*}
This completes the proof of Theorem \ref{thm:BCE}.
When $\sgm{u}$ is trigonometric or rational,
Proposition~\ref{prop:LPBC},~(2) implies that
these functions $\Phi_{BC}(x;y|\delta,\kappa)$
and $\Psi_{BC}(x;y)$ satisfy
the functional equations as stated in Theorem~\ref{thm:BCT}.

\subsection[Dif\/ference operators of Koornwinder type]{Dif\/ference operators of Koornwinder type}\label{section3.4}

In the rest of this section, we conf\/ine ourselves to
the trigonometric and rational $BC$ cases and suppose that
$\sgm{u}$ is a constant multiple of $\sin(\pi u/\omega_1)$ or $u$.
We rewrite our results on kernel functions for these cases,
in terms of dif\/ference operator
\begin{gather*}
D_{x}^{(\mu|\delta,\kappa)}
=
\dsum{i=1}{m}
A_{i}^{+}(x;\mu|\delta,\kappa) \big(T_{x_i}^{\delta}-1\big)
+
\dsum{i=1}{m}
A_{i}^{-}(x;\mu|\delta,\kappa) \big(T_{x_i}^{-\delta}-1\big)
\end{gather*}
{\em of Koornwinder type}.
We remark that this operator has symmetry
\begin{gather*}
D_{x}^{(\mu|-\delta,-\kappa)}=D_{x}^{(-\mu|\delta,\kappa)}
\end{gather*}
with respect to the sign change,
as in the case of $E_x^{(\mu|\delta,\kappa)}$.
We show f\/irst that $D_x^{(\mu|\delta,\kappa)}$
dif\/fers from $E_x^{(\mu|\delta,\kappa)}$ only by
an additive constant in the 0th order term.

\begin{lemma}\qquad
\begin{enumerate}\itemsep=0pt
\item[$(1)$] The constant function $1$ is an eigenfunction of
$E_x^{(\mu|\delta,\kappa)}$:
\begin{gather*}
E_x^{(\mu|\delta,\kappa)}(1)=
C^{(\mu|\delta,\kappa)}+
\dfrac{1}{\sgm{\kappa}}\big(
\sgm{2m\kappa+c^{(\mu|\delta,\kappa)}}-
\sgm{c^{(\mu|\delta,\kappa)}}\big).
\end{gather*}
\item[$(2)$] The two difference operators $D_x^{(\mu|\delta,\kappa)}$
and $E_x^{(\mu|\delta,\kappa)}$ are related as
\begin{gather*}
E_x^{(\mu|\delta,\kappa)}=D_x^{(\mu|\delta,\kappa)}
+
C^{(\mu|\delta,\kappa)}+
\dfrac{1}{\sgm{\kappa}}\big(
\sgm{2m\kappa+c^{(\mu|\delta,\kappa)}}-
\sgm{c^{(\mu|\delta,\kappa)}}
\big).
\end{gather*}
\end{enumerate}
\end{lemma}

\begin{proof}
Since $D_x^{(\mu|\kappa)}=E_x^{(\mu|\delta,\kappa)}-E_x^{(\mu|\delta,\kappa)}(1)$, statement (2) follows from statement (1).
For the proof of (1), we make use of our Theorem \ref{thm:BCT}.
This theorem is valid even in the case where
the dimension $m$ or $n$ reduces
to zero. When $n=0$, Theorem \ref{thm:BCT}, (1) implies
\begin{gather*}
\sgm{\kappa} E_x^{(\mu|\delta,\kappa)}(1)-\sgm{\kappa}
C^{(\nu|\delta,\kappa)}=
\sgm{2m\kappa+c^{(\mu|\delta,\kappa)}},
\end{gather*}
with $\nu_s={\frac 12}(\delta+\kappa)-\mu_s$ ($s=1,\ldots,2\rho$),
since $\Phi_{BC}(x;y)$ in this case is the constant function 1.
Also from the case $m=n=0$ we have
\begin{gather*}
\sgm{\kappa} C^{(\mu|\delta,\kappa)}
-\sgm{\kappa} C^{(\nu|\delta,\kappa)}=\sgm{c^{(\mu|\delta,\kappa)}}.
\end{gather*}
Combining these two formulas we obtain
\begin{gather*}
\sgm{\kappa} E_x^{(\mu|\delta,\kappa)}(1)=
\sgm{\kappa} C^{(\mu|\delta,\kappa)}+
\sgm{2m\kappa+c^{(\mu|\delta,\kappa)}}-
\sgm{c^{(\mu|\delta,\kappa)}}.
\tag*{\qed}
\end{gather*}
\renewcommand{\qed}{}
\end{proof}

Let us rewrite the functional equations in
Proposition \ref{prop:LPBC} in terms of the operator
$D_x^{(\mu|\delta,\kappa)}$. In the notation of Proposition
\ref{prop:LPBC}, (2) we have
\begin{gather}\label{eq:eqmn}
\sgm{\kappa} E_{x}^{(\mu|\delta,\kappa)}\Phi(x;y)
+
\sgm{\lambda} E_{y}^{(\nu|\tau,\lambda)}\Phi(x;y)
=\sgm{2m\kappa+2n\lambda+c^{(\mu|\delta,\kappa)}}\Phi(x;y).
\end{gather}
As special cases where $(m,n)=(m,0),\, (0,n)$ and $(0,0)$, we
have
\begin{gather*}
\sgm{\kappa} E_{x}^{(\mu|\delta,\kappa)}(1)
+
\sgm{\lambda} C^{(\nu|\tau,\lambda)}
 = \sgm{2m\kappa+c^{(\mu|\delta,\kappa)}},\\
\sgm{\kappa} C^{(\mu|\delta,\kappa)}
+
\sgm{\lambda} E_y^{(\nu|\tau,\lambda)}(1)
 = \sgm{2n\lambda+c^{(\mu|\delta,\kappa)}},\\
\sgm{\kappa} C^{(\mu|\delta,\kappa)}
+
\sgm{\lambda} C^{(\nu|\tau,\lambda)}
 = \sgm{c^{(\mu|\delta,\kappa)}},
\end{gather*}
and hence
\begin{gather}
-\sgm{\kappa} E_{x}^{(\mu|\delta,\kappa)}(1)\Phi(x;y)
-
\sgm{\lambda} C^{(\nu|\tau,\lambda)}\Phi(x;y)
 = -\sgm{2m\kappa+c^{(\mu|\delta,\kappa)}}\Phi(x;y),\nonumber\\
-\sgm{\kappa} C^{(\mu|\delta,\kappa)}\Phi(x;y)
-\sgm{\lambda} E_y^{(\nu|\tau,\lambda)}(1)\Phi(x;y)
 = -\sgm{2n\lambda+c^{(\mu|\delta,\kappa)}}\Phi(x;y),\nonumber\\
\sgm{\kappa} C^{(\mu|\delta,\kappa)}\Phi(x;y)
+
\sgm{\lambda} C^{(\nu|\tau,\lambda)}\Phi(x;y)
 = \sgm{c^{(\mu|\delta,\kappa)}}\Phi(x;y),\label{eq:eq00}
\end{gather}
By taking the sum of the four formulas in \eqref{eq:eqmn} and
\eqref{eq:eq00}, we obtain
\begin{gather*}
\sgm{\kappa} D_{x}^{(\mu|\delta,\kappa)}\Phi(x;y)
+
\sgm{\lambda} D_{y}^{(\nu|\tau,\lambda)}\Phi(x;y)
=C \Phi(x;y),
\\
C=
\sgm{2m\kappa + 2n\lambda + c^{(\mu|\delta,\kappa)}}
 - \sgm{2m\kappa + c^{(\mu|\delta,\kappa)}}
 - \sgm{2n\lambda + c^{(\mu|\delta,\kappa)}}
 + \sgm{c^{(\mu|\delta,\kappa)}}.
\end{gather*}
In the rational case, it is clear that $C=0$.
In the trigonometric case, this constant $C$ can be factorized.
In fact, if we choose
$\sgm{u}=2\sqrt{-1}\sin(\pi u/\omega_1)=e(u/2\omega_1)-e(-u/2\omega_1)$,
we have a~simple expression
\begin{gather*}
C=
\sgm{m\kappa}\sgm{n\lambda}\sgm{m\kappa+n\lambda+c^{(\mu|\delta,\kappa)}}.
\end{gather*}
Hence,
under the assumption of Proposition \ref{prop:LPBC}, (2),
we have
\begin{gather*}
\sgm{\kappa}D_x^{(\mu|\delta,\kappa)}\Phi(x;y)
+
\sgm{\lambda}D_y^{(\nu|\tau,\lambda)}\Phi(x;y)
=
\sgm{m\kappa}\sgm{n\lambda}\sgm{m\kappa+n\lambda+c^{(\mu|\delta,\kappa)}}
\Phi(x;y)
\quad(\mbox{resp.} =0),
\end{gather*}
when $\sgm{u}=e(u/2\omega_1)-e(-u/2\omega_1)$
(resp.~when $\sgm{u}=u$).
Applying this to the two cases of \eqref{eq:twocasesBC},
one can easily derive functional equations as in
Theorem~\ref{thm:BCTD}.

\section[Kernel functions for $q$-dif\/ference operators]{Kernel functions for $\boldsymbol{q}$-dif\/ference operators}\label{section4}

In the trigonometric case, it is also important to
consider $q$-dif\/ference operators passing to
multiplicative variables.
Assuming that $\sgm{u}=2\sqrt{-1}\sin(\pi u/\omega_1)
=e(u/2\omega_1)-e(-u/2\omega_1)$,
we def\/ine by $z=e(u/\omega_1)$ the multiplicative
variable associated with the additive variable $u$.
When we write $\sgm{u}=z^{{\frac 12}}-z^{-{\frac 12}}=-z^{-{\frac 12}}(1-z)$,
we regard the square root
$z^{{\frac 12}}$ as a multiplicative notation for $e(u/2\omega_1)$.
We set $q=e(\delta/\omega_1)$ and $t=e(\kappa/\omega_1)$,
assuming that
$\mbox{Im}(\delta/\omega_1)>0$, namely, $|q|<1$.

\subsection{Kernel functions for Macdonald operators}\label{section4.1}

We introduce the multiplicative variables
$z=(z_1,\ldots,z_m)$ and $w=(w_1,\ldots,w_n)$
corresponding to
the additive variables $x=(x_1,\ldots,x_m)$ and
$y=(y_1,\ldots,y_n)$, by
$z_i=e(x_i/\omega_1)$ ($i=1,\ldots,m$) and $w_k=e(y_k/\omega_1)$
($k=1,\ldots,n$),
respectively.

In this convention of the trigonometric case,
the Ruijsenaars dif\/ference operator $D_{x}^{(\delta,\kappa)}$
of type $A$ is a constant multiple of the $q$-dif\/ference operator
of Macdonald:
\begin{gather*}
D_x^{(\delta,\kappa)}=t^{-{\frac 12}(m-1)} \mathcal{D}_z^{(q,t)},
\qquad
\mathcal{D}_z^{(q,t)}
=\dsum{i=1}{m} \dprod{1\le j\le m;\ j\ne i}{}
\dfrac{tz_i-z_j}{z_i-z_j} T_{q,z_i},
\end{gather*}
where $T_{q,z_i}$ denotes the $q$-shift operator with respect to $z_i$:
\begin{gather*}
T_{q,z_i}f(z_1,\ldots,z_m)=f(z_1,\ldots,qz_i,\ldots,z_m)
\qquad(i=1,\ldots,m).
\end{gather*}
In what follows,
we use the gamma function $G_{-}(u|\delta)$ of Section~\ref{section2}, \eqref{eq:tgamma}.
Then our kernel function
$\Phi_{A}(x;y|\delta,\kappa)$ with parameter $v=\kappa$
is expressed as follows
in terms of multiplicative variables:
\begin{gather}
\Phi_{A}(x;y|\delta,\kappa)
=e\big(\tfrac{mn\delta}{2\omega_1}\tbinom{\kappa/\delta}{2}\big)
(z_1\cdots z_m)^{\tfrac{n\kappa}{2\delta}}
(w_1\cdots w_n)^{\tfrac{m\kappa}{2\delta}}
 \Pi(z;w|q,t),\nonumber\\
\Pi(z;w|q,t)=
\dprod{j=1}{m}\dprod{l=1}{n}
\dfrac{(tz_jw_l;q)_\infty}{(z_jw_l;q)_\infty}.\label{eq:KMac}
\end{gather}
The functional equation of Theorem \ref{thm:AT}, (1) thus implies
\begin{gather}\label{eq:CMac}
\mathcal{D}_z^{(q,t)}\Pi(z;w|q,t)-
t^{m-n} \mathcal{D}_w^{(q,t)}\Pi(z;w|q,t)=\dfrac{1-t^{m-n}}{1-t}
\Pi(z;w|q,t),
\end{gather}
which can also be proved by the expansion formula
\begin{gather}\label{eq:CMacExp}
\Pi(z;w|q,t)=\dsum{l(\lambda)\le \min\{m,n\}}{}
b_\lambda(q,t) P_{\lambda}(z|q,t)P_{\lambda}(w|q,t)
\end{gather}
of Cauchy type for Macdonald polynomials~\cite{Ma1995}.
(Formula~\eqref{eq:CMac} already implies
that $\Pi(z;w|q,t)$ has an expansion of this form, apart from
the problem of determining the coef\/f\/icients $b_\lambda(q,t)$.)
On the other hand, kernel function $\Psi_{A}(x;y)$ with parameter
$v=0$ is expressed as
\begin{gather*}
\Psi_{A}(x;y)=(z_1\cdots z_m)^{-\tfrac{n}{2}}
(w_1\cdots w_n)^{-\tfrac{m}{2}}
\dprod{j=1}{m}\dprod{l=1}{n}(z_j-w_l).
\end{gather*}
Then the functional equation
of Theorem \ref{thm:AT}, (2) implies
\begin{gather*}
\big((1-t)\mathcal{D}_{z}^{(q,t)}-(1-q)\mathcal{D}_{w}^{(t,q)}
-(1-t^mq^n)\big)
\dprod{j=1}{m}\dprod{l=1}{n}(z_j-w_l)
=0.
\end{gather*}
This formula corresponds to the dual Cauchy formula
\begin{gather*}
\dprod{j=1}{m}\dprod{l=1}{n}(z_j-w_l)
=\dsum{\lambda\subset(n^m)}{}
(-1)^{|\lambda^\ast|}
P_\lambda(z|q,t) P_{\lambda^\ast}(w|t,q),
\end{gather*}
where $\lambda^\ast=(m-\lambda_n',\ldots,n-\lambda_1')$
is the partition representing the complement of $\lambda$
in the $m\times n$ rectangle.

These kernel functions for Macdonald
operators have been applied to the studies of
raising and lowering operators (Kirillov--Noumi \cite{KN1998, KN1999},
Kajihara--Noumi \cite{KaN2000}) and
integral representation (Mimachi---Noumi \cite{MN1997}, for instance).
We also remark that, in this $A$ type case,
a kernel function of Cauchy type for $q$-Dunkl operators
has been constructed by Mimachi---Noumi \cite{MN1998}.

\subsection{Kernel functions for Koornwinder operators}\label{section4.2}

We now consider the trigonometric $BC$ case.
Instead of the additive parameters
$(\mu_1,\mu_2,\mu_3,\mu_4)$
($\rho=2$), we use the multiplicative parameters
\begin{gather*}
(a,b,c,d)=(e(\mu_1/\omega_1),e(\mu_2/\omega_1),
e(\mu_3/\omega_1),e(\mu_4/\omega_1)).
\end{gather*}
These four parameters are the {\em Askey--Wilson parameters}
$(a,b,c,d)$ for the Koornwinder polynomials
$P_{\lambda}(z;a,b,c,d|q,t)$.

In this trigonometric $BC$ case, the dif\/ference operator
\begin{gather*}
D_{x}^{(\mu|\delta,\kappa)}
=
\dsum{i=1}{m}
A_{i}^{+}(x;\mu|\delta,\kappa) (T_{x_i}^{\delta}-1)
+
\dsum{i=1}{m}
A_{i}^{-}(x;\mu|\delta,\kappa) (T_{x_i}^{-\delta}-1),
\end{gather*}
which we have discussed in Section~\ref{section3.4},
is a constant multiple of Koornwinder's $q$-dif\/ference
operator~\cite{K1992}.
Let us consider the Koornwinder operator
\begin{gather*}
\mathcal{D}_{z}^{(a,b,c,d|q,t)}
=
\dsum{i=1}{m} \mathcal{A}_i^+(z)
 (T_{q,z_i}-1)
+
\dsum{i=1}{m} \mathcal{A}_i^-(z)
 \big(T_{q,z_i}^{-1}-1\big)
\end{gather*}
in the multiplicative variables, where the coef\/f\/icients
$\mathcal{A}_i^{+}(z)=\mathcal{A}_i^{+}(z;a,b,c,d|q,t)$
are given by
\begin{gather*}
\mathcal{A}_i^{+}(z)
=
\dfrac{(1-az_i)(1-bz_i)(1-cz_i)(1-dz_i)}{
(abcdq^{-1})^{\frac 12} t^{m-1} (1-z_i^2)(1-qz_i^2)}
\dprod{1\le j\le m;\ j\ne i}{}
\dfrac{1-tz_iz_j^{\pm1}}{1-z_iz_j^{\pm1}}
\end{gather*}
and $\mathcal{A}_i^{-}(z)=\mathcal{A}_i^{+}(z^{-1})$
for $i=1,\ldots,m$.
Then we have
$D_x^{(\mu|\delta,\kappa)}=-\mathcal{D}_z^{(a,b,c,d|q,t)}$.
Note that this ope\-ra\-tor $\mathcal{D}_z^{(a,b,c,d|q,t)}$ is renormalized
by dividing the one used in~\cite{K1992}
by the factor
$(abcdq^{-1})^{\frac 12} t^{m-1}.\!$
In what follows, we simply suppress the
dependence on the parameters $(a,b,c,d|q,t)$
as $\mathcal{D}_z=\mathcal{D}_z^{(a,b,c,d|q,t)}$,
when we refer to operators or functions associated
with these standard parameters.

In Section~\ref{section2}, we described two types of kernel functions
\eqref{eq:phifunBC} and \eqref{eq:phifunBC2} of
Cauchy type.
Also, depending on the choice of $G(u|\delta)$
we obtain several kernel functions for each type.
From the gamma functions $G_{\mp}(u|\delta)$
of~\eqref{eq:tgamma},
we obtain two kernel functions of type~\eqref{eq:phifunBC};
in the multiplicative variables,
\begin{gather}\label{eq:kerPhi}
\Phi_0(z;w|q,t)
=(z_1\cdots z_m)^{n\beta}
\dprod{j=1}{m}\dprod{l=1}{n}
\dfrac{\big(q^{{\frac 12}}t^{{\frac 12}}z_jw_l^{\pm1};q\big)_\infty}
{\big(q^{\frac 12} t^{-{\frac 12}}z_jw_l^{\pm 1};q\big)_\infty},
\end{gather}
and
\begin{gather*}
\Phi_\infty(z;w|q,t)
=(z_1\cdots z_m)^{-n\beta}
\dprod{j=1}{m}\dprod{l=1}{n}
\dfrac{\big(q^{{\frac 12}}t^{{\frac 12}}z_j^{-1}w_l^{\pm1};q\big)_\infty}
{\big(q^{\frac 12} t^{-{\frac 12}}z_j^{-1}w_l^{\pm 1};q\big)_\infty},
\end{gather*}
respectively,
where we put $\beta=\kappa/\delta$ so that $t=q^\beta$.
Similarly, we obtain two kernel functions of type~\eqref{eq:phifunBC2}
from $G_{\pm}(u|\delta)$:
\begin{gather*}
\Phi_{+}(z;w|q,t)
=e\big(\tfrac{f(x;y)}{\omega_1\delta}\big) \dprod{j=1}{m}\dprod{l=1}{n} \dprod{\epsilon_1,\epsilon_2=\pm}{}
\big(q^{{\frac 12}}t^{{\frac 12}}z_j^{\epsilon_1}w_{l}^{\epsilon_2};q\big)_{\infty},
\\
\Phi_{-}(z;w|q,t)
=e\big( - \tfrac{f(x;y)}{\omega_1\delta}\big) \dprod{j=1}{m}\dprod{l=1}{n}
\dprod{\epsilon_1,\epsilon_2=\pm}{}
\big(q^{{\frac 12}}t^{-{\frac 12}}z_j^{\epsilon_1}w_{l}^{\epsilon_2};q\big)_{\infty}^{-1},
\end{gather*}
where $f(x;y)=n\sum\limits_{j=1}^{m}x_j^2+m\sum\limits_{l=1}^{n} y_l^2
+\tfrac{mn}{4}(\kappa^2-\delta^2)$.
Each of these four functions dif\/fers from the others by multiplicative factors
which are $\delta$-periodic in all the $x$ variables and $y$ variables,
It should be noted, however, that they have dif\/ferent analytic properties.
In the following we denote simply by $\Phi(z;w|q,t)$ one of these
functions.

The kernel function $\Psi(z;w)=\Psi_{BC}(x;y)$
of dual Cauchy type is given by
\begin{gather}\label{eq:kerPsi}
\Psi(z;w)=
\dprod{j=1}{m}\dprod{l=1}{n}
\big(z_j+z_j^{-1}
-w_l-w_l^{-1}\big)
=\dprod{j=1}{m}\dprod{l=1}{n}
(z_j-w_l)\big(1-z_j^{-1}w_l^{-1}\big),
\end{gather}
which is precisely the kernel function introduced by Mimachi~\cite{Mi2001}.

For the passage from additive variables
to multiplicative variables,
we introduce the multiplicative notation
for the function $\sgm{u}$: For $z=e(u/\omega_1)$, we write
$\br{z}=\sgm{u}$. Namely, we~set
\begin{gather*}
\br{z}=z^{{\frac 12}}-z^{-{\frac 12}}=-z^{-{\frac 12}}(1-z), \qquad z=e(u/\omega_1),
\end{gather*}
with the square root $z^{{\frac 12}}$ regarded as the
multiplicative notation for $e(u/2\omega_1)$.
This function $\br{z}$ is a natural object to be used
in the case of $BC$ type,
because of the symmetry $\br{z^{-1}}=-\br{z}$.
In this notation, the coef\/f\/icients $\mathcal{A}_i^{+}(z)$
of the Koornwinder operator $\mathcal{D}_{z}$ are
expressed simply~as
\begin{gather*}
\mathcal{A}_i^{+}(z)
=
\dfrac{\br{az_i}\br{bz_i}\br{cz_i}\br{dz_i}}
{\br{z_i^2}\br{qz_i^2}}
\dprod{1\le j\le m;\, j\ne i}{}
\dfrac{\br{tz_i/z_j}\br{tz_iz_j}}{\br{z_i/z_j}\br{z_iz_j}}
\qquad(i=1,\ldots,m).
\end{gather*}
It should be noted also that our parameter
$c^{(\mu|\delta,\kappa)}
=\sum\limits_{s=1}^{4}\mu_s - (\delta + \kappa) + \omega_1$
passes to multiplicative variables as
\begin{gather*}
\sgm{u+c^{(\mu|\delta,\kappa)}}
=-\left[u+\sum\limits_{s=1}^{4}\mu_s-\delta-\kappa\right]
=-\br{zabcd/qt},\qquad z=e(u/\omega_1),
\end{gather*}
with a minus sign.
Then, Theorem \ref{thm:BCTD} can be restated as follows.
\begin{theorem}\label{thm:CDKoorn}\qquad
\begin{enumerate}\itemsep=0pt
\item[$(1)$] The function $\Phi(z;w|q,t)$ defined as above
satisfies the functional equation
\begin{gather*}
\br{t}
\mathcal{D}_z\Phi(z;w|q,t)-
\br{t}
\widetilde{\mathcal{D}}_w\Phi(z;w|q,t)
=\br{t^m}\br{t^{-n}}\br{abcdq^{-1}t^{m-n-1}}\Phi(z;w|q,t),
\end{gather*}
where $\widetilde{\mathcal{D}}_w$ denotes the Koornwinder operator
in $w$ variables with parameters
$(a,b,c,d)$
replaced by
$(\sqrt{qt}/a,\sqrt{qt}/b,\sqrt{qt}/c,\sqrt{qt}/d)$.

\item[$(2)$] The function $\Psi(z;w)$ defined as~\eqref{eq:kerPsi}
satisfies the functional equation
\begin{gather*}
\br{t}
\mathcal{D}_z\Psi(z;w)
+
\br{q}
\widehat{\mathcal{D}}_w\Psi(z;w)
=\br{t^m}\br{q^n}\br{abcdt^{m-1}q^{n-1}}
\Psi(z;w),
\end{gather*}
where $\widehat{D}_w$ denotes the Koornwinder operator
in $w$ variables with parameters $(a,b,c,d| t,q)$.
\end{enumerate}
\end{theorem}

Statement (2) of Theorem \ref{thm:CDKoorn} recovers
the key lemma of Mimachi \cite[Lemma 3.2]{Mi2001},
from which he established the dual Cauchy formula
\begin{gather*}
\dprod{j=1}{m}\dprod{l=1}{n}
\big(z_j + z_j^{-1} - w_l - w_l^{-1}\big)
=\dsum{\lambda\subset(n^m)}{}
(-1)^{|\lambda^\ast|}
P_{\lambda}(z;a,b,c,d|q,t)
P_{\lambda^\ast}(w;a,b,c,d|t,q)
\end{gather*}
for Koornwinder polynomials,
where the summation is taken over all partitions
$\lambda=(\lambda_1,\ldots,\lambda_m)$ contained
in the $m\times n$ rectangle, and
$\lambda^\ast=(m-\lambda_n',\ldots,m-\lambda_1')$.
By this formula, he also constructed
an integral representation
of Selberg type
for Koornwinder polynomials
attached to rectangles
$(n^m)$ ($n=0,1,2,\ldots$).
We expect that our kernel function $\Phi(z;w|q,t)$
of Cauchy type
could be applied as well to the study of eigenfunctions
of the $q$-dif\/ference operators of Koornwinder.
As a f\/irst step of such applications,
in Section~\ref{section5} we construct explicit
formulas for Koornwinder polynomials attached to
single columns and single rows.

\section{Application to Koornwinder polynomials}\label{section5}

In this section, we apply our results on the kernel functions
for Koornwinder operators to the study of Koornwinder polynomials.
In particular,
we present new explicit formulas for Koornwinder
polynomials attached to single columns and single rows.

To be more precise, we make use of the kernel functions
to express Koornwinder polynomials
$P_{(1^r)}(z;a,b,c,d|q,t)$ ($r=0,1,\ldots,m$)
and
$P_{(l)}(z;a,b,c,d|q,t)$ ($l=0,1,2,\ldots$)
in terms of certain explicitly def\/ined Laurent polynomials
$E_{r}(z;a|t)$ and $H_l(z;a|q,t)$, respectively
(Theorems~\ref{thm:Kcolumn} and~\ref{thm:Krow}).
We remark that these Laurent polynomials
$E_{r}(z;a|t)$ are $H_l(z;a|q,t)$ are in fact
constant multiples
of the $BC_m$ interpolation polynomials
of Okounkov~\cite{O1998} attached
to the partitions $(1^r)$ and $(l)$, respectively.
(This fact will be proved in Appendix~\ref{appendixC}.)
Namely, Theorems~\ref{thm:Kcolumn} and~\ref{thm:Krow}
provide with two special cases of the binomial expansion of the
Koornwinder polynomials in terms of $BC_m$ interpolation
polynomials as is discussed in Okounkov~\cite{O1998}
and Rains~\cite{Ra2005}.

Once we establish the fact that $E_{r}(z;a|t)$ and $H_l(z;a|q,t)$
are interpolation polynomials,
Theorems \ref{thm:Kcolumn} and \ref{thm:Krow} can also be
obtained from Okounkov's binomial formula \cite{O1998},
together with Rains' explicit evaluation of the
binomial coef\/f\/icients for the cases of $(1^r)$ and $(l)$ \cite{Ra2005}\footnote{The authors thank Professor Eric Rains for pointing out
this connection with the interpolation polynomials and
the binomial formula.}.

Before starting the discussion of Koornwinder polynomials,
we introduce a notation
\begin{gather*}
\br{z;w}=\br{zw}\br{z/w}=z+z^{-1}-w-w^{-1}
\end{gather*}
which corresponds to $[u\pm v]=[u+v][u-v]$ in additive variables.
This expression, which appears frequently in the discussion of type
$BC$, deserves a special attention.
Note that
\begin{gather*}
\br{a^{-1};b}=\br{a;b^{-1}}=\br{a;b},\qquad
\br{b;a}=-\br{a;b},\qquad
\br{a;b}+\br{b;c}=\br{a;c},
\end{gather*}
as clearly seen by the def\/inition.
Also, the Riemann relation for $\sgm{u}$ can be written as
\begin{gather*}
\br{z;a}\br{b;c}+\br{z;b}\br{c;a}+\br{z;c}\br{a;b}=0,
\qquad
\dfrac{\br{z;a}}{\br{z;b}}-
\dfrac{\br{w;a}}{\br{w;b}}=
\dfrac{\br{z;w}\br{a;b}}{\br{z;b}\br{w;b}}.
\end{gather*}

\subsection{Koornwinder polynomials}\label{section5.1}

We brief\/ly recall some basic facts about Koornwinder polynomials;
for details, see Stokman \cite{St2004} for example.

Let $\mathbb{K}=\mathbb{Q}\big(a^{\frac 12},b^{\frac 12},c^{\frac 12},d^{\frac 12},q^{\frac 12},t^{\frac 12}\big)$
be the f\/ield of rational functions
in indeterminates,
representing the square roots of the parameters $a$, $b$, $c$, $d$, $q$, $t$,
with coef\/f\/icients in $\mathbb{Q}$, and
$\mathbb{K}[z^{\pm1}]=
\mathbb{K}[z_1^{\pm1},\ldots,z_m^{\pm1}]$ the ring of
Laurent polynomials in $m$ variables
$z=(z_1,\ldots,z_m)$ with coef\/f\/icients in $\mathbb{K}$.
Then the Weyl group
$W=\{\pm1\}^m\rtimes \mathfrak{S}_m$ of type $BC$ (hyperoctahedral group)
acts naturally on $\mathbb{K}[z^{\pm1}]$
through the permutation of indices
for the $z$ variables and the individual inversion of
variables $z_i$ ($i=1,\ldots,m$).
The {\em Koornwinder polynomial}
$P_{\lambda}(z)=P_{\lambda}(z;a,b,c,d|q,t)$
attached to a partition $\lambda=(\lambda_1,\ldots,\lambda_m)$
is then characterized as a unique $W$-invariant Laurent polynomial
in $\mathbb{K}[z^{\pm1}]$ satisfying the following two conditions:
\begin{enumerate}\itemsep=0pt
\item[(1)] $P_{\lambda}(z)$ is expressed in terms
of orbit sums $m_\mu(z)=\sum\limits_{\nu\in W.\mu} z^\nu$
as
\begin{gather*}
P_{\lambda}(z)=m_{\lambda}(z)+\dsum{\mu<\lambda}{} c_{\lambda,\mu}
m_{\mu}(z)\qquad(c_{\lambda,\mu}\in \mathbb{K}),
\end{gather*}
where $\le$ stands for the dominance ordering of
partitions.
\item[(2)] $P_{\lambda}(z)$ is an eigenfunction of Koornwinder's
$q$-dif\/ference operator $\mathcal{D}_z$:
\begin{gather*}
\mathcal{D}_z P_{\lambda}(z)=d_\lambda P_\lambda(z)
\qquad\mbox{for some }\ \ d_\lambda\in \mathbb{K}.
\end{gather*}
\end{enumerate}

These polynomials $P_{\lambda}(z)$, indexed by partitions
$\lambda$, form a $\mathbb{K}$-basis of the ring of $W$-invariants
$\mathbb{K}[z^{\pm1}]^{W}$.
Also, the eigenvalues $d_\lambda$ are given by
\begin{gather*}
d_\lambda =\dsum{i=1}{m} \br{\alpha t^{m-i}q^{\lambda_i};\alpha t^{m-i}}
 =\dsum{i=1}{m}\big(\alpha t^{m-i}q^{\lambda_i}+\alpha^{-1} t^{-m+i}q^{-\lambda_i}
\big)-\dsum{i=1}{m}\big(\alpha t^{m-i}+\alpha^{-1} t^{-m+i}\big),
\end{gather*}
in our notation $\br{z;w}=z + z^{-1} - w - w^{-1}$,
where $\alpha=(abcdq^{-1})^{\frac 12}$. (Note that
$P_{0}(z)=1$, and $d_{0}=0$.)

We give below new explicit formulas for
Koornwinder polynomials $P_{(1^r)}(z)$ attached
to single columns $(1^r)$  ($r=0,1,\ldots,m$), and
$P_{(r)}(z)$ attached to single rows ($r=0,1,2\ldots$).
As we already mentioned, our explicit formulas provide
with the two special cases of Okounkov's binomial expansion
of the Koornwinder polynomials in terms of $BC_m$ interpolation polyno\-mials.
We also remark that, in the cases of type $B$, $C$, $D$, some conjectures
have been proposed by Lassalle~\cite{L2004} on explicit formulas
for Macdonald polynomials attached to single rows.  The relationship
between his conjectures and our approach will be discussed in
a separate paper.

In order to formulate our results, we def\/ine
a set of $W$-invariant Laurent polynomials
$E_r(z;a)$ {\em with reference point $a$}
($r=0,1,\ldots,m$) by
\begin{gather}\label{eq:defpolE}
E_r(z;a|t)=\dsum{1\le i_1<\cdots<i_r\le m}{}
\br{z_{i_1};t^{i_1-1}a}
\br{z_{i_2};t^{i_2-2}a}
\cdots
\br{z_{i_r};t^{i_r-r}a}.
\end{gather}
As we will see below, these Laurent polynomials are $W$-invariant
in spite of their appearance, and they can be considered
as a variation of the orbit sums
$m_{(1^r)}(z)$ ($r=0,1,\ldots,m$)
attached to the fundamental weights.
In fact, these Laurent polynomials $E_r(z;a|t)$ ($r=0,1,\ldots,m$)
are essentially the $BC_m$ interpolation polynomials
of Okounkov attached to single columns $(1^r)$ (for a proof, see Appendix~\ref{appendixC}).
We remark that these polynomials had appeared already in
the work of van Diejen~\cite{vD1996} in relation to the eigenvalues
of his commuting family $q$-dif\/ference operators for this
$BC_m$ case.
They are also used ef\/fectively by a recent work of
Aomoto--Ito~\cite{AI2009} in their study of Jackson integrals
of type $BC$.

\begin{theorem}\label{thm:Kcolumn}
The Koornwinder polynomials $P_{(1^r)}(z;a,b,c,d|q,t)$
attached to columns $(1^r)$ $(r=0,1,\ldots,m)$
are expressed as follows
in terms of $E_{s}(z;a|t)$ $(s=0,1,\ldots,m)$:
\begin{gather}
P_{(1^r)}(z;a,b,c,d|q,t)
=
\dsum{l=0}{r}
\dfrac{
\br{t^{m-r+1}, t^{m-r}ab, t^{m-r}ac, t^{m-r}ad}_{t,l}
}
{\br{t,t^{2(m-r)}abcd}_{t,l}}
E_{r-l}(z;a|t),\label{eq:Kcolumn}
\end{gather}
where $\br{a}_{t,l}=\br{a}\br{ta}\cdots\br{t^{l-1}a}$,
and $\br{a_1,\ldots,a_r}_{t,l}=\br{a_1}_{t,l}\cdots\br{a_r}_{t,l}$.
\end{theorem}
By using $\br{a}_{t,l}=(-1)^lt^{-{\frac 12}\tbinom{l}{2}}a^{-\frac{l}{2}}(a;t)_l$,
formula \eqref{eq:Kcolumn} can be rewritten as follows
in terms of
ordinary $t$-shifted factorials of \cite{GR2004}:
\begin{gather*}
P_{(1^r)}(z;a,b,c,d|q,t)
=
\dsum{l=0}{r}
\dfrac{
(t^{m-r+1}, t^{m-r}ab, t^{m-r}ac, t^{m-r}ad ;t)_l
}
{t^{\tbinom{l}{2}+(m-r)l}a^{l} (t, t^{2(m-r)}abcd ;t)_l}
E_{r-l}(z;a|t).
\end{gather*}

For the description of Koornwinder polynomials attached to single rows,
we introduce a~sequence of $W$-invariant Laurent polynomials
$H_{l}(z;a|q,t)$ ($l=0,1,2,\ldots$) as follows:
\begin{gather}
H_{l}(z;a|q,t)\nonumber\\
=\dsum{\nu_1+\cdots+\nu_m=l}{}
\dfrac{\br{t}_{q,\nu_1}\cdots\br{t}_{q,\nu_m}
}{\br{q}_{q,\nu_1}\cdots\br{q}_{q,\nu_m}}
\br{z_1;a}_{q,\nu_1}
\br{z_2;tq^{\nu_1}a}_{q,\nu_2}
\cdots
\br{z_m;t^{m-1}q^{\nu_1+\cdots+\nu_{m-1}}a}_{q,\nu_m},\label{eq:defpolH}
\end{gather}
where
\begin{gather*}
\br{z;a}_{q,l}=\br{z;a}\br{z;qa}\cdots \br{z;q^{l-1}a}
=(-1)^l q^{-\tbinom{l}{2}}a^{-l}\big(az,az^{-1};q\big)_l.
\end{gather*}
Note that
$\br{z;a}_{q,l}$ is a monic Laurent polynomial in $z$ of
degree $l$.
(The $W$-invariance of $H_l(z;a|q,t)$ will be proved in Lemma \ref{lem:expHw} below.)
These Laurent polynomials $H_{l}(z;a|q,t)$ can be regarded as a $BC_m$ analogue of
the $A_{m-1}$ Macdonald polynomials attached to single rows:
\begin{gather*}
\dfrac{(t;q)_r}{(q;q)_r}P_{(r)}^{A}(z)
=Q_{(r)}^{A}(z)=
\dsum{\nu_1+\cdots+\nu_m=r}{}
\dfrac
{(t;q)_{\nu_1}\cdots(t;q)_{\nu_m}}
{(q;q)_{\nu_1}\cdots(q;q)_{\nu_m}}
z_1^{\nu_1}\cdots z_m^{\nu_m}.
\end{gather*}
Also, they are special cases of $BC_m$ interpolation polynomials attached
to single rows $(l)$ (see Appendix~\ref{appendixC}).

\begin{theorem}\label{thm:Krow}
The Koornwinder polynomials $P_{(r)}(z;a,b,c,d|q,t)$
attached to rows $(r)$ $(r=0,1,2,\ldots)$ are expressed as follows
in terms of $H_{l}(z;a|q,t)$ $(l=0,1,2,\ldots)$:
\begin{gather}
\dfrac{\br{t}_{q,r}}{\br{q}_{q,r}}
P_{(r)}(z;a,b,c,d|q,t)=
\dfrac{
\br{t^m,t^{m-1}ab,t^{m-1}ac,t^{m-1}ad}_{q,r}
}
{
\br{q,t^{2(m-1)}abcdq^{r-1}}_{q,r}
}
\nonumber\\
\phantom{\dfrac{\br{t}_{q,r}}{\br{q}_{q,r}}
P_{(r)}(z;a,b,c,d|q,t)=}{}\times
\dsum{l=0}{r}
\dfrac{
(-1)^l
\br{q^{-r},t^{2(m-1)}abcdq^{r-1}}_{q,l}
}
{
\br{t^m,t^{m-1}ab,t^{m-1}ac,t^{m-1}ad}_{q,l}
}
H_{l}(z;a|q,t).\label{eq:Krow}
\end{gather}
\end{theorem}

Let us denote by $p_r(z;a,b,c,d|q)$ the Koornwinder polynomial
$P_{(r)}(z;a,b,c,d|q,t)$
in the one variable case ($r=0,1,2,\ldots$).
Note that, when $m=1$, $H_l(z;a|q,t)$ reduces to
$\br{t}_{q,l}\br{z;a}_{q,l}/\br{q}_{q,l}$.
Hence Theorem \ref{thm:Krow} for $m=1$ implies
\begin{gather*}
p_r(z;a,b,c,d|q) =
\dfrac{\br{ab, ac, ad}_{q,r}}
{\br{abcdq^{r-1}}_{q,r}}
\dsum{l=0}{r}
\dfrac{(-1)^l\br{q^{-r},abcdq^{r-1}}_{q,l}}
{
\br{q,ab,ac,ad}_{q,l}
}
\br{z;a}_{q,l}\\
\phantom{p_r(z;a,b,c,d|q)}{} =
\dfrac{(ab,ac,ad;q)_r}{a^r(abcdq^{r-1};q)_r}\,
{}_4\phi_{3}
\left[
\begin{matrix}
q^{-r}, abcdq^{r-1}, az, a/z\\
ab,\ ac,\ ad
\end{matrix};q,q
\right],
\end{gather*}
which recovers the well-known ${}_4\phi_{3}$ representation
of the (monic) Askey--Wilson polynomials.

We prove these Theorems \ref{thm:Kcolumn} and \ref{thm:Krow}
in Subsections~\ref{section5.2} and~\ref{section5.3},
by means of the kernel functions of dual Cauchy type,
and of Cauchy type, respectively.

\subsection{Case of a single column}\label{section5.2}

We f\/irst explain some properties of the {\em elementary}
Laurent polynomials
$E_r(z;a|t)$ ($r=0,1$, $\ldots,m$).
\begin{lemma}
The Laurent polynomials $E_r(z;a|t)$ $(r=0,1,\ldots,m)$
are characterized as the expansion coefficients in
\begin{gather*}
\dprod{j=1}{m}\br{w;z_j}=\dsum{r=0}{m}(-1)^r E_r(z;a|t) \br{w;a}_{t,m-r},
\end{gather*}
where $\br{w;a}_{t,l}=\br{w;a}\br{w;ta}\cdots\br{w;t^{l-1}a}$.
In particular,  $E_r(z;a|t)$ is $W$-invariant for each $r=0,1,\ldots,m$.
\end{lemma}

\begin{proof}
Since the uniqueness of expansion in terms of $\br{w;a}_{t,r}$
($r=0,1,\ldots,m$) is obvious, we show the validity of the
expansion formula above. Note that $E_r(z;a|t)$ can be expressed
as
\begin{gather*}
E_r(z;a|t)
=\dsum{|I|=r}{}
\dprod{i\in I}{} \br{z_i;t^{|I^{\rm c}_{<i}|}a},
\qquad I^{\rm c}_{<i}=\pr{j\in \{1,\ldots,m\}\backslash I\ |\ j<i}.
\end{gather*}
Hence we have the recurrence formula
\begin{gather}\label{eq:polErec}
E_{r}(z;a|t)=\br{z_1;a}E_{r-1}(z';a|t)+
E_{r}(z';ta|t),\qquad z'=(z_2,\ldots,z_m).
\end{gather}
By using this recurrence, one can inductively
prove the expansion formula of lemma. In fact,
by $\br{w;z_1}=\br{w;a}-\br{z_1;a}$, we compute
\begin{gather*}
\dprod{j=1}{m}\br{w;z_j}
=\br{w;a}\dprod{i=2}{m}\br{w;z_i}-\br{z_1;a}\dprod{i=2}{m}\br{w;z_i}\\
\hphantom{\dprod{j=1}{m}\br{w;z_j}}{} =\br{w;a}\dsum{r=0}{m-1}(-1)^r E_{r}(z';ta|t)\br{w;ta}_{t,m-1-r}\\
\hphantom{\dprod{j=1}{m}\br{w;z_j}=}{}
-\br{z_1;a}\dsum{s=0}{m-1}(-1)^s E_{s}(z';a|t)\br{w;a}_{t,m-1-s}\\
\hphantom{\dprod{j=1}{m}\br{w;z_j}}{}=\dsum{r=0}{m}(-1)^r \big(E_{r}(z';ta|t)+\br{z_1;a}
E_{r-1}(z';a|t)\big)\br{w;a}_{t,m-r}
\\
\hphantom{\dprod{j=1}{m}\br{w;z_j}}{}=\dsum{r=0}{m}(-1)^r
E_{r}(z;a|t)\br{w;a}_{t,m-r}.\tag*{\qed}
\end{gather*}
\renewcommand{\qed}{}
\end{proof}

We now proceed to the proof of Theorem \ref{thm:Kcolumn}.
Mimachi's dual Cauchy formula for Koornwinder polynomials
can be written as
\begin{gather*}
\dprod{j=1}{m}\dprod{l=1}{n}\br{w_l;z_j}
=\dsum{\lambda\subset(n^m)}{}
(-1)^{|\lambda|}
P_{\lambda}(z;a,b,c,d|q,t)
P_{\lambda^\ast}(w;a,b,c,d|t,q).
\end{gather*}
When $n=1$, this formula implies
\begin{gather}\label{eq:KcolA}
\dprod{j=1}{m} \br{w;z_j}
=\dsum{r=0}{m}
(-1)^{r}
P_{(1^r)}(z|q,t) p_{m-r}(w|t),
\end{gather}
where we have omitted the parameters $(a,b,c,d)$.
Namely, the Koornwinder polynomials attached to single columns are
determined as expansion coef\/f\/icients of
the kernel function for $n=1$ in terms of the monic
Askey--Wilson polynomials
$p_l(w|t)=p_l(w;a,b,c,d|t)$ with base $t$.
On the other hand, we already have the expansion formula
\begin{gather}\label{eq:KcolB}
\dprod{j=1}{m} \br{w;z_j}
=\dsum{l=0}{m}(-1)^l E_{l}(z;a|t) \br{w;a}_{t,m-l}.
\end{gather}
Recalling that
\begin{gather}
p_l(w|t)
 =
\dfrac{\br{ab, ac, ad}_{t,l}}
{\br{abcdt^{l-1}}_{t,l}}
\dsum{r=0}{l}
(-1)^r\dfrac{\br{t^{-l},abcdt^{l-1}}_{t,r}}
{\br{t,ab,ac,ad}_{t,r}
}
\br{w;a}_{t,r}\nonumber\\
\phantom{p_l(w|t)}{}=
\dsum{r=0}{l}
\dfrac{\br{t^{r+1},t^rab,t^rac,t^rad}_{t,l-r}}
{\br{t,abcdt^{l+r-1}}_{t,l-r}}
\br{w;a}_{t,r}
\qquad(l=0,1,2,\ldots)\label{eq:pinw}
\end{gather}
we consider to express $\br{w;a}_{t,l}$ in terms of Askey--Wilson
polynomials $p_r(w|t)$.

\begin{lemma}\label{lem:inv}
For each $l=0,1,\ldots$, one has
\begin{gather}
\br{w;a}_{t,l}
 =\dsum{r=0}{l}
(-1)^{l-r}
\dfrac{
\br{ t^{r+1}, t^rab, t^rac, t^rad}_{t,l-r}
}
{
\br{ t,abcdt^{2r}}_{t,l-r}
}
 p_{r}(w|t)
\qquad(l=0,1,2,\ldots).\label{eq:inv}
\end{gather}
\end{lemma}

We omit the proof of this lemma,
since it can be derived as a special case of the connection formula
for Askey--Wilson polynomials with dif\/ferent parameters
(see \cite{GR2004}). Note that, if we set $d=t^{1-l}/a$ in~\eqref{eq:pinw}, then $p_l(w;a,b,c,t^{1-l}/a|t)=\br{w;a}_{t,l}$.

By substituting \eqref{eq:inv} into \eqref{eq:KcolB}, we obtain
\begin{gather*}
\dprod{j=1}{m} \br{w;z_j}
=\dsum{0\le l\le r\le m}{}(-1)^r
\dfrac{
\br{t^{m-r+1},t^{m-r}ab,t^{m-r}ac,t^{m-r}ad}_{r-l}
}
{
\br{t,abcdt^{2(m-r)}}_{t,r-l}
}
E_{l}(z;a|t)
p_{m-r}(w|t).
\end{gather*}
Comparing this formula with \eqref{eq:KcolA} we obtain
\begin{gather*}
P_{(1^{r})}(z|q,t)
=\dsum{l=0}{r}
\dfrac{
\br{t^{m-r+1},t^{m-r}ab,t^{m-r}ac,t^{m-r}ad}_{r-l}
}
{
\br{t,abcdt^{2(m-r)}}_{t,r-l}
}
E_{l}(z;a|t),
\end{gather*}
as desired.

\subsection{Case of a single row}\label{section5.3}

Recall that the kernel function of Cauchy type
\begin{gather*}
\Phi(z;w|q,t)
=(z_1\cdots z_m)^{n\beta}
\dprod{j=1}{m}\dprod{l=1}{n}
\dfrac{\big(q^{{\frac 12}}t^{{\frac 12}}z_jw_l^{\pm1};q\big)_\infty}
{\big(q^{\frac 12} t^{-{\frac 12}}z_jw_l^{\pm 1};q\big)_\infty},
\end{gather*}
def\/ined in \eqref{eq:kerPhi},
satisf\/ies the dif\/ference equation
\begin{gather*}
\br{t}
\mathcal{D}_z\Phi(z;w|q,t)-
\br{t}
\widetilde{\mathcal{D}}_w\Phi(z;w|q,t)
=\br{t^m}\br{t^{-n}}\br{abcdq^{-1}t^{m-n-1}}\Phi(z;w|q,t),
\end{gather*}
where $\widetilde{\mathcal{D}}_w$ denotes the Koornwinder operator
in $w$ variables with parameters $(a,b,c,d)$
replaced by
$(\sqrt{qt}/a,\sqrt{qt}/b,
\sqrt{qt}/c,\sqrt{qt}/d| q,t)$.
We set hereafter
\begin{gather*}
\widetilde{a}=\sqrt{qt}/a, \qquad
\widetilde{b}=\sqrt{qt}/b, \qquad
\widetilde{c}=\sqrt{qt}/c, \qquad \widetilde{d}=\sqrt{qt}/d.
\end{gather*}
Also,
for any Laurent polynomial $f(z)\in\mathbb{K}[z^{\pm1}]$,
we denote by $\widetilde{f}(z)\in \mathbb{K}[z^{\pm1}]$ the
Laurent polynomial obtained from $f(z)$ by replacing the
parameters $(a,b,c,d)$ with
$(\widetilde{a},\widetilde{b},\widetilde{c},\widetilde{d})$.

Let us consider the special case where
$t=q^{-k}$ ($k=0,1,2,\ldots$). Then
the kernel function $\Phi(z;w|q,q^{-k})$ reduces to
a Laurent polynomial in $(z,w)$:
\begin{gather*}
\Phi\big(z;w|q,q^{-k}\big)
 =(z_1\cdots z_m)^{-kn}
\dprod{j=1}{m}\dprod{l=1}{n}
\big(q^{{\frac 12}(1-k)}z_jw_l^{\pm 1};q\big)_{k}
 =(-1)^{kmn}
\dprod{j=1}{m}\dprod{l=1}{n} \br{w_l;q^{{\frac 12}(1-k)}z_j}_{q,k}.\!
\end{gather*}
Ignoring the sign factor, we set
\begin{gather*}
\Phi_{-k}(z;w)=
\dprod{j=1}{m}\dprod{l=1}{n} \br{w_l;q^{{\frac 12}(1-k)}z_j}_{q,k}.
\end{gather*}
Note that
\begin{gather*}
\br{w;q^{{\frac 12}{(1-k)}}z}_{q,k}=
\br{w;q^{{\frac 12}{(-k+1)}}z}\br{w;q^{{\frac 12}{(-k+2)}}z}
\cdots\br{w;q^{{\frac 12}{(k-1)}}z}
\end{gather*}
is invariant under the inversion $z\to z^{-1}$.
In what follows, we analyze the case where
$n=1$ and $t=q^{-k}$ ($k=0,1,2,\ldots$).
In this case, the kernel function
\begin{gather*}
\Phi_{-k}(z;w)=\dprod{j=1}{m}\br{w;q^{{\frac 12}(1-k)}z_j}_k
\end{gather*}
is a symmetric Laurent polynomial in $w$ of degree $km$.
Also, this kernel function satisf\/ies the functional equation
\begin{gather*}
\mathcal{D}_z\Phi_{-k}(z;w)-
\widetilde{\mathcal{D}}_w\Phi_{-k}(z;w)
=-\br{t^m}\br{abcdq^{-1}t^{m-2}}\Phi_{-k}(z;w).
\end{gather*}
With $\alpha=(abcdq^{-1})^{{\frac 12}}$, this formula can be written
as
\begin{gather}\label{eq:kerrelKrow}
\mathcal{D}_z\Phi_{-k}(z;w)=
\big(\widetilde{\mathcal{D}}_w-\br{\alpha t^{m-1};\alpha t^{-1}}\big)
\Phi_{-k}(z;w).
\end{gather}

Noting that $\Phi_{-k}(z;w)$ is a symmetric Laurent polynomial,
we expand this kernel in terms of the monic Askey--Wilson
polynomials
$\widetilde{p}_l(w|q)=p_l(w;\widetilde{a},\widetilde{b},\widetilde{c},\widetilde{d}|q)$
in $w$ with the {\em twisted} parameters, so that
\begin{gather*}
\Phi_{-k}(z;w)=\dprod{j=1}{m}\br{w;q^{{\frac 12}(1-k)}z_j}_k
=\dsum{l=0}{km} G_l(z) \widetilde{p}_{km-l}(w|q).
\end{gather*}
The Laurent polynomials $G_l(z)$ ($0\le l\le km$) are uniquely
determined by this expansion, and hence $W$-invariant.

\begin{lemma}
When $t=q^{-k}$ $(k=0,1,2,\ldots)$,
the $W$-invariant Laurent polynomials $G_l(z)$ $(0\le l\le km)$
defined as above are eigenfunctions
of $\mathcal{D}_z$:
\begin{gather*}
\mathcal{D}_zG_l(z)=\br{\alpha t^{m-1}q^l;\alpha t^{m-1}}G_l(z)
\qquad (0\le l\le km).
\end{gather*}
\end{lemma}

\begin{proof}
Note that
\begin{gather*}
\mathcal{D}_w p_{km-l}(w|q)=
\br{\alpha q^{km-l};\alpha} p_{km-l}(w|q)
=\br{\alpha t^{-m}q^{-l};\alpha} p_{km-l}(w|q),
\end{gather*}
where $\alpha=(abcdq^{-1})^{{\frac 12}}$.
Hence,
by $\widetilde{\alpha}=tq^{{\frac 12}}(abcd)^{-{\frac 12}}=t/\alpha$, we
obtain
\begin{gather*}
\widetilde{\mathcal{D}}_w
\widetilde{p}_{km-l}(w|q)=
\br{t^{1-m}q^{-l}/\alpha;t/\alpha}
\widetilde{p}_{km-l}(w|q)
=\br{\alpha t^{m-1}q^l;\alpha t^{-1}}\widetilde{p}_{km-l}(w|q).
\end{gather*}
In view of \eqref{eq:kerrelKrow}, we compute
\begin{gather*}
\big(\widetilde{\mathcal{D}}_w -\br{\alpha t^{m-1}; \alpha t^{-1}}\big)
\widetilde{p}_{km-l}(w|q)
=\big(\br{\alpha t^{m-1}q^l;\alpha t^{-1}}
-\br{\alpha t^{m-1}; \alpha t^{-1}}\big)
\widetilde{p}_{km-l}(w|q)
\\
\phantom{\big(\widetilde{\mathcal{D}}_w -\br{\alpha t^{m-1}; \alpha t^{-1}}\big)
\widetilde{p}_{km-l}(w|q)}  {} =\br{\alpha t^{m-1}q^l;\alpha t^{m-1}}
\widetilde{p}_{km-l}(w|q).
\end{gather*}
Hence, \eqref{eq:kerrelKrow} implies
\begin{gather*}
\dsum{l=0}{km} \mathcal{D}_zG_l(z) \widetilde{p}_{km-l}(w|q)
=
\dsum{l=0}{km} G_l(z)
\big(\widetilde{\mathcal{D}}_w -\br{\alpha t^{m-1}; \alpha t^{-1}}\big)\widetilde{p}_{km-l}(w|q)\\
\phantom{\dsum{l=0}{km} \mathcal{D}_zG_l(z) \widetilde{p}_{km-l}(w|q)}{} =
\dsum{l=0}{km}
\br{\alpha t^{m-1}q^l;\alpha t^{m-1}}
G_l(z)
\widetilde{p}_{km-l}(w|q).
\end{gather*}
Namely,
\begin{gather*}
\mathcal{D}_zG_l(z)=\br{\alpha t^{m-1}q^l;\alpha t^{m-1}}G_l(z)
\qquad (l=0,1,\ldots,km). \tag*{\qed}
\end{gather*}
\renewcommand{\qed}{}
\end{proof}

As we have seen above,
each $G_l(z)$ ($l=0,1,\ldots,km$) is an eigenfunction of
$\mathcal{D}_z$ with precisely the same eigenvalue as
the one for the Koornwinder polynomial
$P_{(l)}(z)$ attached to the single row of length $l$.
At this moment, however, we cannot conclude this
$G_{l}(z)$ {\em is} indeed a constant multiple of the
Koornwinder polynomial $P_{(l)}(z)$ specialized to the case
$t=q^{-k}$. This is because dif\/ferent partitions $\lambda$
may give the same eigenvalue
\begin{gather*}
d_\lambda=\sum\limits_{i=1}^{m}\br{\alpha t^{m-i}q^{\lambda_i};\alpha t^{m-i}}
\end{gather*}
under this specialization. This point will be discussed later
after we determine an explicit formula for $G_l(z)$.

Since we already know the relationship between
the Askey--Wilson polynomials $p_{l}(w|q)$ and the
Laurent polynomials $\br{w;a}_{q,l}$, we consider
to expand the kernel function $\Phi_{-k}(z;w)$ in
terms of $\br{w;a}_{q,l}$.

\begin{lemma}\label{lem:expHw}\qquad
\begin{enumerate}\itemsep=0pt
\item[$(1)$]
When $t=q^{-k}$ $(k=0,1,2,\ldots)$, the kernel function $\Phi_{-k}(z;w)$
has the following expansion formula$:$
\begin{gather}\label{eq:expHw}
\Phi_{-k}(z;w)=\dprod{j=1}{m}\br{w;q^{{\frac 12}(1-k)}z_j}_{q,k}
=\dsum{l=0}{km} H_l(z) \br{w;\sqrt{qt}/a}_{q,km-l},
\end{gather}
where $H_l(z)$ stands for $H_l(z;a|q,t)$
defined as \eqref{eq:defpolH} with $t=q^{-k}$ for $l=0,1,\ldots,km$.
\item[$(2)$]
The Laurent polynomials $H_l(z;a|q,t)$ $(l=0,1,2,\ldots)$
are $W$-invariant.
\end{enumerate}
\end{lemma}

\begin{proof}
Statement (2) follows from the expansion formula
\eqref{eq:expHw} of statement (1).
Since \linebreak $\Phi_{-k}(z;w)$ is $W$-invariant in the $z$ variables,
formula \eqref{eq:expHw} implies that $H_l(z;a|q,q^{-k})$ is $W$-invariant
for $k\ge l/m$.  Hence we see that
$H_l(z;a|q,t)$ itself is $W$-invariant as a Laurent polynomial
in $\mathbb{K}[z^{\pm}]$ for each $l=0,1,2,\ldots$.

In the following proof of statement (1),
we omit the base $q$, and write $\br{w;a}_{l}=\br{w;a}_{q,l}$.
We f\/irst present a connection formula for
the Laurent polynomials $\br{w;a}_{l}$ and $\br{w;b}_{l}$
with dif\/ferent reference point $a$, $b$:
\begin{gather*}
\br{w;b}_{l}
=\dsum{r=0}{l}(-1)^r\sbinom{l}{r}
 \br{q^{l-r}ab,b/a}_{r} \br{w;a}_{l-r},\qquad
\sbinom{l}{r}=(-1)^r\dfrac{\br{q^{-l}}_r}{\br{q}_r},
\end{gather*}
which is equivalent to the $q$-Saalsch\"utz sum \cite{GR2004}
\begin{gather*}
\dfrac{(bw,b/w;q)_l}{(ba,b/a;q)_l}
={}_3\phi_{2}
\left[
\begin{matrix}
q^{-l}, aw, a/w\\
ab, q^{1-l}a/b
\end{matrix}; q,q
\right].
\end{gather*}
Let us rewrite the formula above in the form
\begin{gather*}
\br{w;b}_{l}
=
\dsum{r=0}{l}(-1)^r\sbinom{l}{r}
 \br{q^{{\frac 12}(l-1)}b;q^{{\frac 12}(1-l)}/a}_{r} \br{w;a}_{l-r}.
\end{gather*}
Hence,
\begin{gather*}
\br{w;q^{{\frac 12}(1-\mu)}z}_{\mu}
  =
\dsum{\nu=0}{\mu}(-1)^\nu\sbinom{\mu}{\nu}
 \br{z;q^{{\frac 12}(1-\mu)}/a}_{\nu} \br{w;a}_{\mu-\nu}.
\end{gather*}
This implies
\begin{gather*}
\br{w;q^{{\frac 12}(1-\mu_1)}z_1}_{\mu_1}
\br{w;q^{{\frac 12}(1-\mu_2)}z_2}_{\mu_2}
\\
\qquad{}=
\dsum{\nu_1=0}{\mu_1}(-1)^{\nu_1}\sbinom{\mu_1}{\nu_1}
\br{z_1;q^{{\frac 12}(1-\mu_1)}/a}_{\nu_1} \br{w;a}_{\mu_1-\nu_1}
\br{w;q^{{\frac 12}(1-\mu_2)}z_2}_{\mu_2}.
\end{gather*}
Then in each term, we expand
$\br{w;q^{{\frac 12}(1-\mu_2)}z_2}_{\mu_2}$
in terms of $\br{w;q^{\mu_1-\nu_1}a}_{\mu_2-\nu_2}$
($0\le\nu_2\le\mu_2$)
to get
\begin{gather*}
\br{w;q^{{\frac 12}(1-\mu_1)}z_1}_{\mu_1}
\br{w;q^{{\frac 12}(1-\mu_2)}z_2}_{\mu_2}\\
\qquad{} =
\dsum{\nu_1=0}{\mu_1}\;
\dsum{\nu_2=0}{\mu_2}
(-1)^{\nu_1+\nu_2}\sbinom{\mu_1}{\nu_1}\sbinom{\mu_2}{\nu_2}
 \br{z_1;q^{{\frac 12}(1-\mu_1)}/a}_{\nu_1}
\\
\qquad \qquad{}\times
 \br{z_2;q^{{\frac 12}(1-\mu_2)-(\mu_1-\nu_1)}/a}_{\nu_2}
\br{w;a}_{\mu_1+\mu_2-\nu_1-\nu_2}.
\end{gather*}
By repeating this procedure, we f\/inally obtain
\begin{gather*}
\dprod{j=1}{m}
\br{w;q^{{\frac 12}(1-\mu_j)}z_j}_{\mu_j}
=
\dsum{\nu}{}(-1)^{|\nu|}
\dprod{j=1}{m}
\sbinom{\mu_j}{\nu_j}
\br{z_j;q^{{\frac 12}(1-\mu_j)-\sum\limits_{i<j}(\mu_i-\nu_i)}/a}_{\nu_j}
\br{w;a}_{|\mu|-|\nu|}
\end{gather*}
for any $\mu=(\mu_1,\ldots,\mu_m)$,
where the sum is taken over all multi-indices $\nu=(\nu_1,\ldots,\nu_m)$
such that $\nu_i\le\mu_i$ ($i=1,\ldots,m$).
As a special case of this formula where $\mu_1=\cdots=\mu_m=k$,
$|\mu|=km$,
we get the expansion formula
\begin{gather*}
\Phi_{-k}(z;w)=\dprod{j=1}{m}
\br{w;q^{{\frac 12}(1-k)}z_j}_{k}
=
\dsum{l=0}{km} \widetilde{H}_{l}(z) \br{w;a}_{km-l}
\end{gather*}
for $\Phi_{-k}(z;w)$.
Here the coef\/f\/icients are determined as
\begin{gather*}
\widetilde{H}_{l}(z)
 =
\dsum{|\nu|=l}{}(-1)^{|\nu|}
\dprod{j=1}{m}
\sbinom{k}{\nu_j} \br{z_j;q^{{\frac 12}(1-k)-k(j-1)+\sum\limits_{i<j}\nu_i}/a}_{\nu_j}
\\
 \phantom{\widetilde{H}_{l}(z)}{} =
\dsum{\nu_1+\cdots+\nu_m=l}{}\;
\dprod{j=1}{m}
\dfrac{\br{t}_{\nu_j}}{\br{q}_{\nu_j}}
\br{z_j;t^{j-1}q^{\sum\limits_{i<j}\nu_i}\sqrt{qt}/a}_{\nu_j}
\end{gather*}
with $t=q^{-k}$.
Replacing the parameters $a$ by
$\widetilde{a}=\sqrt{qt}/a$,
we obtain
\begin{gather*}
\Phi_{-k}(z;w)=\dprod{j=1}{m}
\br{w;q^{{\frac 12}(1-k)}z_j}_{k}
=
\dsum{l=0}{km} H_{l}(z) \br{w;\sqrt{qt}/a}_{km-l},
\end{gather*}
where
\begin{gather*}
H_{l}(z)
=\dsum{\nu_1+\cdots+\nu_m=l}{}
\dprod{j=1}{m}
\dfrac{\br{t}_{\nu_j}}{\br{q}_{\nu_j}}
\br{z_j;t^{j-1}q^{\sum\limits_{i<j}\nu_i}a}_{\nu_j}
\qquad(t=q^{-k}).
\tag*{\qed}
\end{gather*}
\renewcommand{\qed}{}
\end{proof}

We now have two expansions of the kernel function $\Phi_{-k}(z;w)$:
\begin{gather}\label{eq:expGH}
\Phi_{-k}(z;w)=
\dsum{r=0}{km} G_r(z) \widetilde{p}_{km-r}(w|q)
=
\dsum{l=0}{km} H_l(z) \br{w;\sqrt{qt}/a}_{q,km-l}.
\end{gather}
Also, from Lemma \ref{lem:inv} we see
\begin{gather*}
\br{w;a}_{q,l}
 =\dsum{r=0}{l}
(-1)^{l-r}
\dfrac{
\br{q^{r+1},q^rab,q^rac,q^rad}_{q,l-r}
}
{
\br{q,abcdq^{2r}}_{q,l-r}
}
p_{r}(w|q),
\end{gather*}
and hence,
\begin{gather*}
\br{w;\sqrt{qt}/a}_{q,l}
 =\dsum{r=0}{l}
(-1)^{l-r}
\dfrac{
\br{q^{r+1},tq^{r+1}/ab,tq^{r+1}/ac,tq^{r+1}/ad}_{q,l-r}
}
{
\br{q,t^2q^{2(r+1)}/abcd}_{q,l-r}
}
\widetilde{p}_{r}(w|q)
\end{gather*}
for $l=0,1,2,\ldots$.
Substituting this into \eqref{eq:expGH}, we obtain
an expression of $G_r(z)$ in term of~$H_l(z)$ as follows:
\begin{gather*}
G_r(z)
=\dsum{l=0}{r}
(-1)^{r-l}
\dfrac{
\br{q^{km-r+1},tq^{km-r+1}/ab,tq^{km-r+1}/ac,tq^{km-r+1}/ad}_{q,r-l}
}
{
\br{q,t^2q^{2(km-r+1)}/abcd}_{q,r-l}
}
H_l(z)
\\
\phantom{G_r(z)}{} =\dsum{l=0}{r}
(-1)^{r-l}
\dfrac{
\br{q^{1-r}t^{-m},q^{1-r}t^{1-m}/ab,
q^{1-r}t^{1-m}/ac,q^{1-r}t^{1-m}/ad}_{q,r-l}
}
{
\br{q,q^{2(1-r)}t^{2(1-m)}/abcd}_{q,r-l}
}
H_l(z)
\\
\phantom{G_r(z)}{} =
\dfrac{\br{t^{m},t^{m-1}ab,
t^{m-1}ac,t^{m-1}ad}_{q,r}}
{\br{q,t^{2(m-1)}abcd q^{r-1}}_{q,r}}
\dsum{l=0}{r}
\dfrac
{(-1)^{l}\br{q^{-r},t^{2(m-1)}abcdq^{r-1}}_{q,l}}
{\br{t^{m},t^{m-1}ab,
t^{m-1}ac,t^{m-1}ad}_{q,l}}
H_l(z).
\end{gather*}

From the expression obtained above, it is clear that
\begin{gather*}
H_l(z)=\dfrac{\br{t}_{q,l}}{\br{q}_{q,l}}m_{(l)}(z)+\mbox{terms
lower than $(l)$ with respect to $\le$},
\end{gather*}
and
\begin{gather*}
G_r(z)=\dfrac{\br{t}_{q,r}}{\br{q}_{q,r}}m_{(r)}(z)+\mbox{terms
lower than $(r)$ with respect to $\le$}.
\end{gather*}
Note that $\br{t}_{q,r}=\br{q^{-k}}_{q,r}\ne 0$ for $0\le r\le k$.
Also, we already know that each $G_r(z)$ ($r=0,1,\ldots,km$)
satisf\/ies the dif\/ference
equation
\begin{gather*}
\mathcal{D}_{z}G_{r}(z)=\br{\alpha t^{m-1}q^r;\alpha t^{m-1}} G_{r}(z).
\end{gather*}
Suppose in general that a partition
$\lambda=(\lambda_1,\ldots,\lambda_m)$ satisf\/ies the condition
$\lambda_i-\lambda_{i+1}\le k$ ($i=1,\ldots,m-1$).
Since
\begin{gather*}
\lambda_1+k\le \lambda_2+2k\le \cdots\le\lambda_m+km
\end{gather*}
in this case, it turns out that the eigenvalue
$d_{\mu}$ for any partition $\mu<\lambda$ is distinct
from $d_{\lambda}$
when the square roots of $a$, $b$, $c$, $d$, $q$ are regarded as indeterminates.
For such a partition~$\lambda$, the Koornwinder polynomial
$P_{\lambda}(z)=P_{\lambda}(z;a,b,c,d|q,t)$ can be specialized to
$t=q^{-k}$, and any eigenfunction having the nontrivial
leading term $m_{\lambda}(z)$ must be a constant multiple
of $P_{\lambda}(z;a,b,c,d|q,q^{-k})$.
This implies that, for each $r$ with $0\le r\le k$,
$G_r(z)$ is a constant multiple of $P_{(r)}(z)$
specialized to $t=q^{-k}$:
\begin{gather*}
G_r(z)=\dfrac{\br{t}_{q,r}}{\br{q}_{q,r}}P_{(r)}(z)\big|_{t=q^{-k}}
\qquad (0\le r\le k).
\end{gather*}

For each $r=0,1,2,\ldots$, consider the
Laurent polynomial in $\mathbb{K}[z^{\pm1}]$ def\/ined
by right-hand side of the explicit formula \eqref{eq:Krow}
of Theorem \ref{thm:Krow}.
Then the both sides of \eqref{eq:Krow} are regular at $t=q^{-k}$
($k\ge r$), and they coincide with each other
for $t=q^{-k}$ ($k=r,r+1,\ldots$). Hence the both sides must
be identical as rational functions in $t^{{\frac 12}}$.
This completes the proof of Theorem \ref{thm:Krow}.



\appendix


\section[Remarks on higher order dif\/ference operators]{Remarks on higher order dif\/ference operators}\label{appendixA}

In the case of type $A$, an explicit commuting family of higher order
dif\/ference operators, denoted below by $D_{r,x}^{(\delta,\kappa)}$
($r=1,\ldots,m$), including $D_{x}^{(\delta,\kappa)}$ as a f\/irst member,
has been constructed by Ruijsenaars \cite{Ru1987}.  He also proved in \cite{Ru2009}
that, when $m=n$, the kernel function of Cauchy type, corresponding to our
$\Phi_{A}(x;y|\delta,\kappa)$,
intertwines the whole commuting families of dif\/ference operators in $x$ variable
and $y$ variables.  (For the comparison of our $\Phi_{A}(x;y|\delta,\kappa)$ with
the type~$A$ kernel functions of Ruijsenaars \cite{Ru2009}, see Appendix~\ref{appendixB.1}.)

Fix any nonzero entire function $\sgm{x}$ satisfying the Riemann relation
as in Section~\ref{section2}.
We consider a sequence of
dif\/ference operators $D_{r,x}^{(\delta,\kappa)}$ ($r=1,\ldots,m$) def\/ined by
\begin{gather}\label{eq:ADr}
D_{r,x}^{(\delta,\kappa)}
=\dsum{I\subset\{1,\ldots,m\};\,|I|=r}{} \;
\dprod{i\in I,\,  j\notin I}{}
\dfrac{\sgm{x_i-x_j+\kappa}}{\sgm{x_i-x_j}}\dprod{i\in I}{}T_{x_i}^{\delta}.
\end{gather}
Then, from the result of \cite{Ru1987} and its degenerate cases,
it follows that these operators
$D_{r,x}^{(\delta,\kappa)}$ ($r=1,\ldots,m$) commute with each other.
In this setting the same kernel function
\begin{gather*}
\Phi_{A}(x;y|\delta,\kappa)=\dprod{j,l=1}{m}\dfrac{G(x_j+y_l+v-\kappa|\delta)}
{G(x_j+y_l+v|\delta)}
\end{gather*}
as in \eqref{eq:phifun} for the case $m=n$ satisf\/ies the dif\/ference equation
\begin{gather}\label{eq:AKRr}
D_{r,x}^{(\delta,\kappa)}
\Phi_{A}(x;y|\delta,\kappa)=D_{r,y}^{(\delta,\kappa)}\Phi_{A}(x;y|\delta,\kappa)
\end{gather}
for all $r=1,\ldots,m$.
This functional equation is in fact equivalent to
\begin{gather*}
\dsum{|I|=r}{}\;
\dprod{i\in I;\, j\notin I}{}
\dfrac{\sgm{x_i-x_j+\kappa}}{\sgm{x_i-x_j}}\dprod{i\in I; \, 1\le l\le m}{}
\dfrac{\sgm{x_i+y_l+v-\kappa}}{\sgm{x_i+y_l+v}}\\
\qquad {}=
\dsum{|K|=r}{}\;
\dprod{k\in K;\,l\notin K}{}
\dfrac{\sgm{y_k-y_l+\kappa}}{\sgm{y_k-y_l}}\dprod{k\in K; \, 1\le j\le m}{}
\dfrac{\sgm{y_k+x_j+v-\kappa}}{\sgm{y_k+x_j+v}},
\end{gather*}
which is precisely the key identity of Kajihara--Noumi \cite[Theorem~1.3]{KN2003},
that was derived from the determinantal formula of Frobenius.

In the trigonometric case, it is convenient to consider the generating function
\begin{gather*}
\mathcal{D}^{(q,t)}_z(u)=\dsum{r=0}{m}(-u)^r \mathcal{D}_{r,z}^{(q,t)},\qquad
\mathcal{D}^{(q,t)}_{r,z}=t^{\binom{r}{2}}\dsum{|I|=r}{}\;
\dprod{i\in I; j\notin I}{}
\dfrac{tz_i-z_j}{z_i-z_j}\dprod{i\in I}{}T_{q,z_i}
\end{gather*}
of the Macdonald $q$-dif\/ference operators, passing to the multiplicative variables.
Then from the eigenfunction expansion \eqref{eq:CMacExp}, it follows
that the kernel function
$\Pi(z;w|q,t)$ of \eqref{eq:KMac} for the variables
$z=(z_1,\ldots,z_m)$ and $w=(w_1,\ldots,w_n)$ satisf\/ies the functional equation
\begin{gather*}
(t^mu;t)_\infty \mathcal{D}_{z}^{(q,t)}(t^n u) \Pi(z;w|q,t) =
(t^nu;t)_\infty \mathcal{D}_{w}^{(q,t)}(t^m u) \Pi(z;w|q,t).
\end{gather*}
It would be an important problem to f\/ind an elliptic extension of this formula
for the case $m\ne n$.

As for type $BC$,
a commuting family of higher order dif\/ference operators
for $E_{x}^{(\mu|\delta,\kappa)}$ has been constructed
explicitly by van Diejen \cite{vD1994, vD1996} in the trigonometric case, and
inductively by Komori--Hikami \cite{KH1997} in the elliptic case.
We expect that our kernel function $\Phi_{BC}(x;y|\delta,\kappa)$ should
intertwine the whole commuting families of higher order dif\/ference operators in
$x$ variables and~$y$~variables (at least under the balancing condition in the
elliptic case), similarly to the $A$ type case.

\section{Comparison with \cite{Ru2009}}\label{appendixB}

As we already mentioned,
in his series of works \cite{Ru2009} Ruijsenaars has constructed kernel
function for elliptic dif\/ference operators of type $(A_{m-1},A_{m-1})$ and $(BC_m,BC_m)$
(see also \cite{K1992,L2004}).
In this section, we clarify how our dif\/ference operators and kernel functions
in the elliptic case are related to those of Ruijsenaars.

In what follows, we conf\/ine ourselves to the elliptic case, and
set
\begin{gather}\label{eq:A1}
p=e(\omega_2/\omega_1), \quad q=e(\delta/\omega_1);\qquad
\sgm{x}=-z^{-{\frac 12}}\theta(z;p),\quad z=e(x/\omega_1)
\end{gather}
as in Section~\ref{section2.1}, (2) Elliptic case,
assuming that
$\mbox{Im}(\omega_2/\omega_1)>0$, $\mbox{Im}(\delta/\omega_1)>0$.
We remark that the quasi-periodicity of the function $\sgm{x}$ is described as
\begin{gather*}
\sgm{x+\omega_r}=\epsilon_r e(\eta_r(x+\tfrac{\omega_r}{2})) \sgm{x}
\qquad(r=0,1,2,3),
\end{gather*}
where $\omega_0=0$, $\omega_3=-\omega_1-\omega_2$, and
\begin{gather*}
\epsilon_0=1,\qquad
\epsilon_1=\epsilon_2=
\epsilon_3=-1;\qquad
\eta_0=\eta_1=0,\qquad\eta_2=-\tfrac{1}{\omega_1}, \qquad \eta_3=\tfrac{1}{\omega_1}.
\end{gather*}
(We use below the index 0 instead of 4 for $\omega_r$, $\epsilon_r$ and $\eta_r$.)
Since
$\prod\limits_{s=1}^{3}\sgm{{\frac 12} \omega_s}=-2 e(-\tfrac{\omega_2}{2\omega_1})=
-2p^{-{\frac 12}}$ in this case, we have the duplication formula
\begin{gather*}
\sgm{2x}=2e\big(\tfrac{\omega_2}{2\omega_1}\big) \dprod{r=0}{3}\big[x-\tfrac{\omega_r}{2}\big].
\end{gather*}
In def\/ining kernel functions, we use the gamma function
\begin{gather*}
G(x|\delta)=G_{+}(x|\delta)=e\big(\tfrac{\delta}{2\omega_1}\tbinom{x/\delta}{2}\big)
\Gamma(pz;p,q)
\end{gather*}
of \eqref{eq:egammapm}, associated with $\sgm{x}$.

In the works of Ruijsenaars \cite{Ru2009}, the two periods $\omega_1$, $\omega_2$ and
the scaling constant $\delta$ are
parametrized as
\begin{gather*}
\omega_1=``\dfrac{\pi}{r}",\qquad \omega_2=``ia_+",\qquad \delta=``ia_-"
\end{gather*}
in view of the symmetry between $\omega_2$ and $\delta$ (or $p$ and $q$).
In terms of the $R$-function
\begin{gather*}
\cR(x)=``R(r,a_{+};x)"=\theta\big(p^{{\frac 12}}z;p\big)
\end{gather*}
def\/ined in \cite[I, (1.21)]{Ru2009}, our $\sgm{x}$ is expressed as
\begin{gather*}
[x]=-e\big(-\tfrac{x}{2\omega_1}\big) \cR\big(x-\tfrac{\omega_2}{2}\big).
\end{gather*}
Note that $\cR(-x)=\cR(x)$.
The elliptic gamma function
\begin{gather*}
\cG(x)=``G(r,a_+,a_-;x)"=\Gamma\big(p^{\frac 12} q^{\frac 12} z;p,q\big)
\end{gather*}
of Ruijsenaars \cite[I, (1.19)]{Ru2009}
is related to our $G(x|\delta)$ by the formula
\begin{gather*}
G(x|\delta)=
e\big(\tfrac{\delta}{2\omega_1}\tbinom{x/\delta}{2}\big)
\cG\big(x+{\tfrac 12}(\omega_2-\delta)\big).
\end{gather*}
The function $\cG(x)$ is symmetric with respect to $\omega_2$ and $\delta$
($p$ and $q$), and satisf\/ies the functional equation
\begin{gather*}
\cG\big(x+{\tfrac 12}\delta\big)=\cR(x) \cG\big(x-{\tfrac 12}\delta\big), \qquad \cG(-x)=\cG(x)^{-1}.
\end{gather*}

\subsection[Case of type $A$]{Case of type $\boldsymbol{A}$}\label{appendixB.1}

As in \eqref{eq:ADr},
we consider the commuting family of dif\/ference operators $D_{r,x}^{(\delta,\kappa)}$
($r=1,\ldots,m$) of type $A$ in the variables $x=(x_1,\ldots,x_m)$.
These operators are in fact identical to the dif\/ference operator
$A_{r,+}(-x)$ of \cite[I, (2.1)]{Ru2009}, up to multiplication by constants:
\begin{gather*}
e\big(\tfrac{r(m-r)\kappa}{2\omega_1}\big)
D_{r,x}^{(\delta,\kappa)}
 =
\dsum{|I|=r}{}\;
\dprod{i\in I;\,j\notin I}{}
\dfrac{\cR(x_i-x_j+\kappa-\frac{\omega_2}{2})}{\cR(x_i-x_j-\tfrac{\omega_2}{2}))}
\exp\left(\delta\sum_{i\in I}{}\partial_{x_i}\right)
=``A_{r,+}(-x)",
\end{gather*}
under the identif\/ication of the parameter $``\mu"=\kappa$;
we denote below this operator by $\mathcal{A}_{r,-x}=``A_{r,+}(-x)"$.
In the multiplicative variables
$z_i=e(x_i/\omega_1)$ $(i=1,\ldots,m)$, $t=e(\kappa/\omega_1)$,
the same operator is expressed as
\begin{gather*}
\mathcal{A}_{r,-x}=t^{{\frac 12} r(m-r)}
D_{r,x}^{(\delta,\kappa)}
=
\dsum{|I|=r}{}\;
\dprod{i\in I;\,j\notin I}{}
\dfrac{\theta(tz_i/z_j;p)}{\theta(z_i/z_j;p)} \dprod{i\in I}{}T_{q,z_i}.
\end{gather*}
The kernel function $\Phi_{A}(x;y|\delta,\kappa)$ of \eqref{eq:phifun}
for the $m=n$ case can be expressed in terms of~$\cG(x)$~as
\begin{gather*}
\Phi_{A}(x;y|\delta,\kappa)
=
e\left(-\tfrac{m\kappa}{2\omega_1\delta}\left(\sum_{j=1}^{m}x_j+\sum_{l=1}^{m}y_l\right)
+\tfrac{m^2}{4\omega_1\delta}\kappa(\kappa+\delta-2v)\right)
\mathcal{S}_A(-x;y),
\\
\mathcal{S}_A(-x;y) =
\dprod{j,l=1}{m}
\dfrac{\cG(-x_j-y_l-v+{\frac 12}(\delta-\omega_2))}
{\cG(-x_j-y_l-v+\kappa+{\frac 12}(\delta-\omega_2))}\\
 \phantom{\mathcal{S}_A(-x;y)}{} =
\dprod{j,l=1}{m}
\dfrac{\Gamma(pcz_jw_l/t;p,q)}
{\Gamma(pcz_jw_l;p,q)}\qquad(c=e(v/\omega_1)).
\end{gather*}
The second factor $\mathcal{S}_{A}(-x;y)$
of $\Phi_{A}(x;y|\delta,\kappa)$ coincides with
$``\mathcal{S}_\xi(-x,y)" $ of \cite[I, (2.6)]{Ru2009},
under the identif\/ication of parameters
$``\mu"=\kappa$ and $``\xi"=-v+{\frac 12}(\delta+\kappa-\omega_2)$.
Since the f\/irst factor of $\Phi_{A}(x;y|\delta,\kappa)$ is an eigenfunction
of $T_{x_i}^{\delta}$ and $T_{y_k}^{\delta}$
with an equal eigenvalue $e(-m\kappa/2\omega_1)$,
from~\eqref{eq:AKRr} we obtain
$\mathcal{A}_{r,-x}\mathcal{S}_{A}(-x;y)=\mathcal{A}_{r,-y}\mathcal{S}_{A}(-x;y)$,
namely
\begin{gather*}
\mathcal{A}_{r,x}\mathcal{S}_{A}(x;y)=\mathcal{A}_{r,-y}\mathcal{S}_{A}(x;y)
\qquad(r=1,\ldots,m).
\end{gather*}
This gives formula~(2.5) in \cite[I]{Ru2009} with $``\delta=+"$.
The statement for $``\delta=-"$ follows from the symmetry of $\mathcal{S}_{A}(x;y)$
between~$\omega_2$ and~$\delta$ (or~$p$ and~$q$).

\subsection[Case of type $BC$]{Case of type $\boldsymbol{BC}$}\label{appendixB.2}

We now consider the dif\/ference operator $E_x^{(\mu|\delta,\kappa)}$ of type $BC_m$,
def\/ined by \eqref{eq:defE}, \eqref{eq:defE1}, \eqref{eq:defE2}
with the function $\sgm{x}$ of \eqref{eq:A1}.
For the variables $x=(x_1,\ldots,x_m)$
and the parameters $\mu=(\mu_1,\ldots,\mu_8)$, $\kappa$,
we will use the multiplicative expressions together by setting
\begin{gather*}
z_i=e(x_i/\omega_1)\quad (i=1,\ldots,m),\qquad
u_k=e(\mu_k/\omega_1)\quad (k=1,\ldots,8),\qquad
t=e(\kappa/\omega_1).
\end{gather*}

We f\/irst remark that, by \eqref{eq:defEE},
the dif\/ference operator $E_{x}^{(\mu|\delta,\kappa)}$ can be rewritten
in the form
\begin{gather}
t^{m-1}p^{-1}q^{-\frac{3}{2}}(u_1\cdots u_8)^{{\frac 12}}
E_{x}^{(\mu|\delta,\kappa)}
\nonumber\\
\qquad {} =
\dsum{1\le i\le m;\, \epsilon=\pm}{}
\big(qz_i^{\epsilon 2}\big)^{-1}\cA_{i}^{\epsilon}(x;\mu|\delta,\kappa) T_{x_i}^{\epsilon \delta}
+\dsum{r=0}{3} \cA_r^{0}(x;\mu|\delta,\kappa).\label{eq:AppBC1}
\end{gather}
Here the coef\/f\/icients
$\cA_{i}^{\pm}(x;\mu|\delta,\kappa)$ ($i=1,\ldots,m$),
$\cA_{r}^{0}(x;\mu|\delta,\kappa)$ ($r=0,1,2,3$)
are obtained as follows by modifying the corresponding coef\/f\/icients of
$E_{x}^{(\mu|\delta,\kappa)}$:
\begin{gather*}
\cA_i^{+}(x;\mu|\delta,\kappa)
 =
\dfrac{\prod\limits_{k=1}^{8}\cR\big(x_i + \mu_k - \frac{\omega_2}{2}\big)}
{\cR\big(2x_i - \frac{\omega_2}{2}\big) \cR\big(2x_i + \delta - \frac{\omega_2}{2}\big)}
\dprod{1\le j\le m;\,j\ne i}{}
\dfrac{\cR\big(x_i\pm x_j + \kappa - \frac{\omega_2}{2}\big)}
{\cR\big(x_i\pm x_j - \frac{\omega_2}{2})}\\
\phantom{\cA_i^{+}(x;\mu|\delta,\kappa)}{} =
\dfrac{\prod\limits_{k=1}^{8}\theta(u_kz_i;p)}{\theta(z_i^2;p)\theta(qz_i^2;p)}
\dprod{1\le j\le m;\,j\ne i}{}
\dfrac{\theta(tz_iz_j^{\pm1};p)}{\theta(z_iz_j^{\pm1};p)},
\\
\cA_{i}^{-}(x;\mu|\delta,\kappa)
 =
\cA_{i}^{+}(-x;\mu|\delta,\kappa)\qquad(i=1,\ldots,m),
\\
\cA_r^{0}(x;\mu|\delta,\kappa)
 =
C_r\dfrac{\prod\limits_{k=1}^{8}\cR\big(\mu_k + \tfrac{\omega_r}{2} - {\frac 12}(\omega_2 + \delta)\big)}
{2 \cR\big(\kappa - \frac{\omega_2}{2}\big) \cR\big(\kappa - \delta - \frac{\omega_2}{2}\big)}
\dprod{j=1}{m}
\dfrac{\cR\big(\pm x_j+ \kappa + \frac{\omega_r}{2} - {\frac 12}(\omega_2 + \delta)\big)}
{\cR\big(\pm x_j + \frac{\omega_r}{2} - {\frac 12}(\omega_2 + \delta)\big)}
\\
\phantom{\cA_r^{0}(x;\mu|\delta,\kappa)}{} =
C_r\dfrac{\prod\limits_{k=1}^{8}\theta\big(u_kq^{-\frac 12}c_r;p\big)}
{2 \theta(t;p) \theta(tq^{-1};p)}
\dprod{j=1}{m}
\dfrac{\theta\big(tq^{-\frac 12}c_rz_j^{\pm1};p\big)}
{\theta\big(q^{-\frac 12}c_rz_j^{\pm1};p\big)}\qquad(r=0,1,2,3)
\end{gather*}
with constants
\begin{gather*}
c_r=e(\omega_r/2\omega_1),\qquad
C_r=e\left(-\tfrac{2\omega_r}{\omega_1}-
\eta_r\left(\omega_r-2\delta+m\kappa+{\tfrac 12}\sum_{k=1}^{8}\mu_k\right)\right)\qquad(r=0,1,2,3).
\end{gather*}
In the multiplicative expression, these constants are expressed as
\begin{gather*}
c_0=1, \qquad  c_1=-1, \qquad c_2=p^{\frac 12}, \qquad c_3=-p^{-\frac 12};
\\
C_0=C_1=1,\qquad
C_2=t^m p^{-1}q^{-2}(u_1\cdots u_8)^{{\frac 12}},\qquad
C_3=t^{-m} p^{3}q^{2}(u_1\cdots u_8)^{-{\frac 12}}.
\end{gather*}
In view of formula \eqref{eq:AppBC1}, we take a nonzero meromorphic
function $P(x)$ satisfying the system of dif\/ference equations
\begin{gather*}
T_{x_i}^{\delta}(P(x))=qz_i^2 P(x)\qquad(i=1,\ldots,m).
\end{gather*}
The simplest choice for such a function is given by
$P(x)=e\big(\frac{1}{\omega_1\delta}\sum\limits_{j=1}^{m}x_j^2\big)$.
Then formula \eqref{eq:AppBC1} implies that
the dif\/ference operator $E_{x}^{(\mu|\delta,\kappa)}$ is expressed
in the form
\begin{gather}\label{eq:relE1}
E_{x}^{(\mu|\delta,\kappa)}
=
t^{-m+1}pq^{\frac{3}{2}}(u_1\cdots u_8)^{-{\frac 12}}
P(x) \cE_{x}^{(\mu|\delta,\kappa)} P(x)^{-1},
\end{gather}
as a constant multiple of the conjugation of a dif\/ference operator
\begin{gather*}
\cE_{x}^{(\mu|\delta,\kappa)}=
\dsum{1\le i\le m;\, \epsilon=\pm}{}
\cA_{i}^{\epsilon}(x;\mu|\delta,\kappa) T_{x_i}^{\epsilon \delta}
+\dsum{r=0}{3} \cA_r^{0}(x;\mu|\delta,\kappa)
\end{gather*}
by $P(x)$.

We also remark that the operator $\cE_{x}^{(\mu|\delta,\kappa)}$
has the following property with respect to shifting~$\kappa$ to~$\kappa-\omega_2$
($t$ to $t/p$):
\begin{gather*}
\cE_{x}^{(\mu|\delta,\kappa)}
=\big(p^2qt^{-2}\big)^{m-1} P(x)^{m-1} \cE_{x}^{(\mu|\delta,\kappa-\omega_2)}
P(x)^{-m+1},\qquad t=e(\kappa/\omega_1).
\end{gather*}
Together with \eqref{eq:relE1}, this also implies
\begin{gather*}
E_{x}^{(\mu|\delta,\kappa)}
=
t^{-3(m-1)}p^{2m-1}q^{m+{\frac 12}}(u_1\cdots u_8)^{-{\frac 12}}
 P(x)^m \cE_{x}^{(\mu|\delta,\kappa-\omega_2)}P(x)^{-m}.
\end{gather*}
Setting $\kappa=\lambda+\omega_2$ ($t=sp$, $s=e(\lambda/\omega_1)$),
we rewrite this formula as
\begin{gather}\label{eq:relE3}
E_{x}^{(\mu|\delta,\lambda+\omega_2)}
=
s^{-3(m-1)}p^{-m+2}q^{m+{\frac 12}}(u_1\cdots u_8)^{-{\frac 12}}
 P(x)^m \cE_{x}^{(\mu|\delta,\lambda)}P(x)^{-m}.
\end{gather}

The dif\/ference operator $\cE_{x}^{(\mu|\delta,\lambda)}$ appearing
in \eqref{eq:relE3} is
essentially the same as the operator $``A_{+}(h,\mu;x)"$ of type
$BC_m$ def\/ined by Ruijsenaars \cite[I, (4.1)--(4.3)]{Ru2009} (see also
(3.1)--(3.9)); the dif\/ference between the two is only by an additive
constant. In fact, we have
\begin{gather}\label{eq:compAE}
``A_{+}(h,\mu;x)"=\cE_{x}^{(\mu|\delta,\lambda)}
-\dsum{r=0}{3}\cA_{r}^{0}(\xi_r,\ldots,\xi_r;\mu|\delta,\lambda),
\end{gather}
where $\xi_0=\xi_2=\omega_1/2$, $\xi_1=\xi_3=0$, under the
identif\/ication of the parameters
\begin{gather*}
``h_k"=\mu_{k}-{\tfrac 12}(\omega_2+\delta)\quad (k=0,1,\ldots,7),
\qquad ``\mu"=\lambda,
\end{gather*}
with $\mu_0=\mu_8$.

By Theorem \ref{thm:BCE} we know that the kernel function of \eqref{eq:phifunBC2}
\begin{gather*}
\Phi_{BC}(x;y|\delta,\kappa)
=
\dprod{j=1}{m}\dprod{l=1}{n}\,\dprod{\epsilon_1,\epsilon_2=\pm}{}
G\big(\epsilon_1x_j+\epsilon_2y_l+{\tfrac 12}(\delta-\kappa)|\delta\big)
\end{gather*}
in the variables $x=(x_1,\ldots,x_m)$ and $y=(y_1,\ldots,y_n)$
satisf\/ies the dif\/ference equation
\begin{gather*}
E_{x}^{(\mu|\delta,\kappa)}\Phi_{BC}(x;y|\delta,\kappa)
=
E_{y}^{(\widetilde{\mu|}\delta,\kappa)}\Phi_{BC}(x;y|\delta,\kappa),\\
\widetilde{\mu}_k={\tfrac 12}(\delta+\kappa)-\mu_k\qquad(k=1,\ldots,8)
\end{gather*}
under the balancing condition
\begin{gather}\label{eq:balance}
(m - n - 1)\kappa - \delta+{\tfrac 12}\sum\limits_{k=1}^{8}\mu_k=0
\qquad(t^{m-n-1}q^{-1}(u_1\cdots u_8)^{{\frac 12}}=1).
\end{gather}
With $G(x|\delta)=G_{+}(x|\delta)$, \eqref{eq:egammapm},
this kernel function is given explicitly by
\begin{gather*}
\Phi_{BC}(x;y|\delta,\kappa)
 =
e\big(\tfrac{f(x;y)}{\omega_1\delta}\big)
\dprod{j=1}{m}
\dprod{l=1}{n}
\,
\dprod{\epsilon_1,\epsilon_2=\pm}{}
{\Gamma\big(pq^{{\frac 12}}t^{-{\frac 12}}z_j^{\epsilon_1}w_l^{\epsilon_2};p,q\big)}
\end{gather*}
where $w_l=e(y_l/\omega_1)$ ($l=1,\ldots,n$), and
$f(x;y)=n\sum\limits_{j=1}^{m}x_j^2+m\sum\limits_{l=1}^{n}{y_l}^2
+\tfrac{mn}{4}(\kappa^2-\delta^2)$.
We express
$\Phi_{BC}(x;y|\delta,\kappa)$ in the form
\begin{gather*}
\Phi_{BC}(x;y|\delta,\kappa)
=\mbox{const}\cdot P(x)^n P(y)^m \mathcal{S}_{BC}(x;y),
\end{gather*}
where
\begin{gather*}
P(x)=e\left(\tfrac{1}{\omega_1\delta}\sum\limits_{j=1}^{m}x_j^2\right),\qquad
P(y)=e\left(\tfrac{1}{\omega_1\delta}\sum\limits_{l=1}^{n}y_l^2\right),
\\
\mathcal{S}_{BC}(x;y)
=
\dprod{j=1}{m}
\dprod{l=1}{n}\,
\dprod{\epsilon_1,\epsilon_2=\pm}{}
{\Gamma\big(pq^{{\frac 12}}t^{-{\frac 12}}z_j^{\epsilon_1}w_l^{\epsilon_2};p,q\big)}.
\end{gather*}
It should be noted that, by the substitution $\kappa=\lambda+\omega_2$
($t=sp$, $s=e(\lambda/\omega_1)$),
the function
\begin{gather*}
\mathcal{S}_{BC}(x;y) =
\dprod{j=1}{m}
\dprod{l=1}{n}\,
\dprod{\epsilon_1,\epsilon_2=\pm}{}
{\Gamma\big(p^{{\frac 12}}q^{{\frac 12}}s^{-{\frac 12}}z_j^{\epsilon_1}w_l^{\epsilon_2};p,q\big)}  =
\dprod{j=1}{m}
\dprod{l=1}{n}\,
\dprod{\epsilon_1,\epsilon_2=\pm}{}
\cG\big(\epsilon_1x_j+\epsilon_2 y_l-{\tfrac 12}\lambda\big)
\end{gather*}
becomes manifestly invariant under the action of the Weyl group of type
$BC_m$ (resp. $BC_n$) on the variables $x=(x_1,\ldots,x_m)$ (resp.\
$y=(y_1,\ldots,y_n)$), and symmetric with respect to
$\omega_2$ and $\delta$ ($p$ and $q$).
This function $\mathcal{S}_{BC}(x;y)$ can be thought
of as the $(BC_m,BC_n)$-version of Ruijsenaars' kernel function
$``\mathcal{S}(x;y)"$ in \cite{Ru2009}.

We now have the functional equation
\begin{gather*}
E_{x}^{(\mu|\delta,\kappa)}P(x)^n P(y)^m \mathcal{S}_{BC}(x;y)
=
E_{y}^{(\widetilde{\mu}|\delta,\kappa)}P(x)^n P(y)^m \mathcal{S}_{BC}(x;y)
\end{gather*}
under the balancing condition \eqref{eq:balance}.
Hence,
by applying \eqref{eq:relE3} with $\kappa=\lambda+\omega_2$, we obtain
\begin{gather*}
\mbox{const}\cdot P(x)^m
\mathcal{E}_{x}^{(\mu|\delta,\lambda)}
P(x)^{n-m} P(y)^{m}\mathcal{S}_{BC}(x;y)
\\
\qquad =
\mbox{const}\cdot P(y)^n
\mathcal{E}_{y}^{(\widetilde{\mu}|\delta,\lambda)}
P(x)^{n} P(y)^{m-n} \mathcal{S}_{BC}(x;y),
\end{gather*}
namely,
\begin{gather*}
\mbox{const}\cdot
\mathcal{E}_{x}^{(\mu|\delta,\lambda)}
P(x)^{n-m} P(y)^{m-n} \mathcal{S}_{BC}(x;y)\\
\qquad =
\mbox{const}\cdot \mathcal{E}_{y}^{(\widetilde{\mu}|\delta,\lambda)}
P(x)^{n-m} P(y)^{m-n} \mathcal{S}_{BC}(x;y).
\end{gather*}
The constants in front of the both sides simplif\/ies by~\eqref{eq:balance},
to imply the functional equation
\begin{gather*}
(pq/s)^{m}
\mathcal{E}_{x}^{(\mu|\delta,\lambda)}
P(x)^{n-m} P(y)^{m-n} \mathcal{S}_{BC}(x;y)
\nonumber\\
\qquad{}=
(pq/s)^{n}
\mathcal{E}_{y}^{(\widetilde{\mu}|\delta,\lambda)}
P(x)^{n-m} P(y)^{m-n} \mathcal{S}_{BC}(x;y),
\end{gather*}
under the balancing condition
\begin{gather*}
{\tfrac 12}\sum\limits_{k=1}^{8}\mu_k
=\delta - (m - n - 1)(\lambda+\omega_2)\qquad
((u_1\cdots u_8)^{{\frac 12}}=q(ps)^{-m+n+1}).
\end{gather*}
In particular, when $m=n$ we have
\begin{gather}\label{eq:Rmm}
\mathcal{E}_{x}^{(\mu|\delta,\lambda)}
\mathcal{S}_{BC}(x;y)
=
\mathcal{E}_{y}^{(\widetilde{\mu}|\delta,\lambda)}
\mathcal{S}_{BC}(x;y),
\end{gather}
under the balancing condition
\begin{gather}\label{eq:Rmmbc}
{\tfrac 12}\sum\limits_{k=1}^{8}\mu_k=\lambda+\omega_2+\delta
\qquad((u_1\cdots u_8)^{{\frac 12}}=spq).
\end{gather}

The functional equation \eqref{eq:Rmm} is precisely
the formula (4.24) of \cite[I, Proposition 4.1]{Ru2009},
for $``\delta=+"$ with the parameter
\begin{gather*}
``h_k"=\mu_k-{\tfrac 12}(\omega_2+\delta)\quad (k=0,1,\ldots,7),\qquad ``\mu"=\lambda.
\end{gather*}
The balancing condition \eqref{eq:Rmmbc} corresponds to
condition~(4.26),
$``\mu=2ia+{\frac 12}\sum\limits_{k=0}^{7}h_k"$.
Also the parameters $\widetilde{\mu}_k={\frac 12}(\lambda+\omega_2+\delta)-\mu_k$
($k=1,\ldots,8$) for the $y$ variables are consistent with
$``(-J_{R}h)_k=\frac{\mu}{2}-ia-h_k"$ ($k=0,1,\ldots,7$).
(The constant $``\sigma_{+}(h)"$ on the right side of \cite[I, (4.24)]{Ru2009},
arises as the dif\/ference of constant terms of two operators, as
indicated in~\eqref{eq:compAE}.)

\section[$E_r(z;a|t)$ and $H_l(z;a|q,t)$ as interpolation polynomials]{$\boldsymbol{E_r(z;a|t)}$ and $\boldsymbol{H_l(z;a|q,t)}$ as interpolation polynomials}\label{appendixC}

In Section~\ref{section5} we presented some explicit expansion formulas for Koornwinder polynomials
attached to single columns and to single rows, in terms of invariant Laurent polynomials
$E_r(z;a|t)$ and $H_l(z;a|q,t)$, respectively.
These Laurent polynomials
$E_r(z;a|t)$ and $H_l(z;a|q,t)$ are essentially the same objects as
the $BC_m$-interpolation polynomials of
Okounkov \cite{O1998} attached to single columns and single rows.
Our Theorems \ref{thm:Kcolumn} and \ref{thm:Krow} provides
with explicit expressions for the corresponding special cases
of the binomial expansion of Koornwinder polynomials in terms
of $BC_m$-interpolation polynomials
as discussed in Okounkov \cite{O1998} and Rains \cite{Ra2005}.

In the notation of Section~\ref{section5},
the $BC_m$ interpolation polynomial $P^\ast_{\lambda}(z;q,t,a)$
attached a~partition $\lambda=(\lambda_1,\ldots,\lambda_m)$ \cite{O1998}
is characterized uniquely, up to multiplication by a constant, as
a Laurent polynomial in $\mathbb{K}[z^{\pm}]$ of degree $|\lambda|$
satisfying the following conditions:
\begin{enumerate}\itemsep=0pt
\item[(0)] $P_\lambda^\ast(z;q,t,a)$ is $W$-invariant in
the shifted variables
$zt^{\delta}a = (z_1t^{m-1}a,z_2t^{m-2}a,\ldots,z_ma)$, where $\delta=(m-1,m-2.\ldots,0)$.
\item[(1)] $P_\lambda^\ast(q^{\mu};q,t,a)=0$ for any partition $\mu=(\mu_1,\ldots,\mu_m)$
such that $\mu\not\supset\lambda$.

\item[(2)] $P_\lambda^\ast(q^{\lambda};q,t,a)\ne 0$.
\end{enumerate}

In this section we show that the Laurent polynomials
$E_{r}(zt^{\delta}a;a|t)$
and $H_{l}(zt^{\delta}a;a|q,t)$
coincide, up to constant multiplication,
with the interpolation polynomials
$P_{(1^r)}^\ast(z;q,t,a)$ and $P_{(l)}^\ast(z;q,t,a)$, respectively.
We prove that these polynomials actually have the interpolation properties
as mentioned above.

\begin{proposition}
For each $r=0,1,\ldots,m$, let $E_r(z;a|t)$ be the $W$-invariant Laurent polynomial
defined in~\eqref{eq:defpolE}.
\begin{enumerate}\itemsep=0pt
\item[$(0)$] $E_r(z_1,\ldots,z_{r},t^{m-r-1}a,\ldots,a;a|t)=\br{z_1;t^{m-r}a}\cdots\br{z_r;t^{m-r}a}.$

\item[$(1)$] For any partition $\mu$ with $l(\mu)<r$, i.e., $\mu\not\supset(1^r)$,
$E_r(q^{\mu}t^{\delta}a;a|t)=0$.

\item[$(2)$]
$E_r(qt^{m-1}a,\ldots,qt^{m-r}a,t^{m-r-1}a,\ldots,a;a|t)
=(-1)^r\br{t^{m-r}a;qt^{m-r}a}_{t,r}.$
\end{enumerate}
\end{proposition}

\begin{proof}
Since $E_r(z;a|t)$ is symmetric in $z=(z_1,\ldots,z_m)$,
from the recurrence relation \eqref{eq:polErec} we have
\begin{gather*}
E_r(z;a|t)=E_{r}(z_1,\ldots,z_{m-1};ta|t)+E_{r-1}(z_1,\ldots,z_{m-1};a|t)\br{z_m;a},
\end{gather*}
and hence by setting $z_m=a$,
\begin{gather*}
E_r(z_1,\ldots,z_{m-1},a;a|t)=E_{r}(z_1,\ldots,z_{m-1};ta|t).
\end{gather*}
This implies
\begin{gather*}
E_r\big(z_1,\ldots,z_{r},t^{m-r-1}a,\ldots,a;a|t\big) =
E_{r}\big(z_1,\ldots,z_{r},t^{m-r-1}a,\ldots,ta;ta|t\big)
 =\cdots\\
\phantom{E_r\big(z_1,\ldots,z_{r},t^{m-r-1}a,\ldots,a;a|t\big)}{} =E_{r}\big(z_1,\ldots,z_{r};t^{m-r}a|t\big)
 =\br{z_1;t^{m-r}a}\cdots\br{z_r;t^{m-r}a}
\end{gather*}
for $r=0,1,\ldots,m$, which proves (0).
In particular, we have
\begin{gather*}
E_r\big(z_1,\ldots,z_{r-1},t^{m-r}a,t^{m-r-1}a,\ldots,a;a|t\big)=0,
\end{gather*}
which implies
$E_r(q^{\mu}t^{\delta}a;a|t)=0$
for any partition such that $l(\mu)<r$, namely, $\mu\not\supset (1^r)$.
Statement (0) also implies
\begin{gather*}
E_r\big(qt^{m-1}a,\ldots,qt^{m-r}a,t^{m-r-1}a,\ldots,a;a|t\big)=(-1)^r\br{t^{m-r}a;qt^{m-r}a}_{t,r}.
\tag*{\qed}
\end{gather*}
\renewcommand{\qed}{}
\end{proof}

From this proposition we see that
$E_r(zt^{\delta}a;a|t)$ is a constant multiple
of the interpolation polynomial $P_{(1^r)}^\ast(z;q,t,a)$
for each $r=0,1,\ldots,m$.

\begin{proposition} For $l=0,1,2,\ldots$,
let $H_l(z;a|q,t)$ be the $W$-invariant
Laurent polynomial defined in~\eqref{eq:defpolH}.
\begin{enumerate}\itemsep=0pt
\item[$(1)$]
For any partition $\mu=(\mu_1,\ldots,\mu_m)$ with $\mu_1<l$, i.e., $\mu\not\supset (l)$,
$H_l(q^{\mu}t^\delta a;a|q,t)=0$.
\item[$(2)$] $H_l(q^lt^{m-1}a,t^{m-2}a,\ldots,a;a|q,t)
=\br{t}_{q,l}\br{q^l t^{2(m-1)}a^2}_{q,l}$.
\end{enumerate}
\end{proposition}

\begin{proof}
By the def\/inition \eqref{eq:defpolH}, the Laurent polynomials
$H_l(z;a|q,t)$ ($l=0,1,2\ldots$) satisfy the recurrence formula
\begin{gather*}
H_l(z;a|q,t)=
\dsum{r=0}{l}\dfrac{\br{t}_{q,l-r}}{\br{q}_{q,l-r}}
 H_{r}(z_1,\ldots,z_{m-1};a|q,t) \br{z_m;q^rt^{m-1}a}_{q,l-r}.
\end{gather*}
Since these Laurent polynomials are symmetric in $z=(z_1,\ldots,z_m)$
(Lemma~\ref{lem:expHw}),
we also have
\begin{gather*}
H_l(z;a|q,t)=
\dsum{r=0}{l}\dfrac{\br{t}_{q,l-r}}{\br{q}_{q,l-r}}
\br{z_1;t^{m-1}q^ra}_{q,l-r}
 H_{r}(z_2,\ldots,z_{m};a|q,t).
\end{gather*}
We prove proposition by the induction on $m$.
For an arbitrary partition $\mu=(\mu_1,\ldots,\mu_m)$,
we specialize this recurrence formula to $z=q^{\mu}t^{\delta}a$:
\begin{gather}
H_l(q^{\mu}t^{\delta}a;a|q,t)
=\dsum{r=0}{l}\dfrac{\br{t}_{q,l-r}}{\br{q}_{q,l-r}}
\br{q^{\mu_1}t^{m-1}a;q^rt^{m-1}a}_{q,l-r}
H_{r}(q^{\mu_2}t^{m-2}a,\ldots,q^{\mu_m}a;a|q,t).\!\!\!\!\label{eq:dainyu}
\end{gather}
By the induction hypothesis we have
\begin{gather*}
H_{r}(q^{\mu_2}t^{m-2}a,\ldots,q^{\mu_m}a;a|q,t)=0
\end{gather*}
for $r>\mu_1$.  For $r\le \mu_1$, we see
\begin{gather*}
\br{q^{\mu_1}t^{m-1}a;q^rt^{m-1}a}_{q,l-r}
=(-1)^{l-r}\br{q^{-\mu_1+r}}_{q,l-r}\br{q^{\mu_1+r}t^{2(m-1)}a^2}_{q,l-r}.
\end{gather*}
Since $-\mu_1+r\le 0$, $\br{q^{-\mu_1+r}}_{q,l-r}=0$ if $l-r>\mu_1-r$, i.e.,
$\mu_1<l$.
This means that all the terms on the right side of \eqref{eq:dainyu}
vanish when $\mu_1<l$, and hence,
$H_l(q^{\mu}t^{\delta}a;a|q,t)=0$.
When $\mu=(l)$, the only nontrivial term arises from $r=0$:
\begin{gather*}
H_l(q^{l}t^{m-1}a,t^{m-2}a,\ldots,a;a|q,t)
=\dfrac{\br{t}_{q,l}}{\br{q}_{q,l}}
\br{q^{l}t^{m-1}a;t^{m-1}a}_{q,l}
=\br{t}_{q,l}\br{q^lt^{2(m-1)}a^2}_{q,l}.
\tag*{\qed}
\end{gather*}
\renewcommand{\qed}{}
\end{proof}

This proposition implies that $H_{l}(zt^{\delta}a;a|a,t)$
is a constant multiple of the interpolation polynomial
$P_{(l)}^{\ast}(z;q,t,a)$ for each $l=0,1,2,\ldots$.

\subsection*{Acknowledgements}

The authors express their sincere gratitude to Professor Simon Ruijsenaars
for his valuable comments concerning the comparison of results in
the present paper with those in his series of works, and to Professor Eric Rains
for suggesting the connection of the explicit formulas for
the Koornwinder polynomials in this paper with
the binomial expansion formula in terms of interpolation polynomials.
Also, they are grateful to the anonymous referees who kindly of\/fered helpful
suggestions for improving the manuscript.

\LastPageEnding

\end{document}